\documentclass[letter, 11pt]{article}
\usepackage{amsmath}

\usepackage[utf8]{inputenc}
\usepackage{array}
\usepackage{cite}
\usepackage{amsfonts}
\usepackage{amssymb}

\usepackage{amsthm}

\usepackage{graphicx} 
\usepackage{epstopdf}
\usepackage{float}
\usepackage{balance}

\usepackage[fancythm,fancybb]{jphmacros2e}
\usepackage{amsfonts}
\usepackage{mathrsfs}

\usepackage{amsthm}

\usepackage{relsize}
\usepackage{capt-of}
\usepackage{mathtools}

\def\frqed{\ifhmode\nobreak\hbox to5pt{\hfil}\nobreak%
	\hskip 0pt plus1fill\nobreak\fi\quad\qedsymbol\renewcommand{\qed}{}} 

\allowdisplaybreaks

\begin{document}



\theoremstyle{plain}
\newtheorem{theorem}{Theorem}[section]
\newtheorem{assumption}[theorem]{Assumption}
\newtheorem{corollary}[theorem]{Corollary}
\newtheorem{lemma}[theorem]{Lemma}
\newtheorem{proposition}[theorem]{Proposition}
\newtheorem{definition}[theorem]{Definition}
\newtheorem{remark}[theorem]{Remark}

\title{Tracking Control by the Newton-Raphson Method with Output Prediction and Controller Speedup}
\author{Y. Wardi, C. Seatzu, J. Cort\'es, M. Egerstedt, S. Shivam, and I. Buckley
\thanks{Wardi, Egerstedt, Shivam and Buckley are with the School of Electrical and Computer Engineering, Georgia Institute of Technology, Atlanta, GA, $30332$, USA.
email: \{ywardi, magnus\}@ece.gatech.edu, \{sshivam6, ihbuckl\}@gatech.edu.\newline
$\-$ $\-$ $\-$  Seatzu is with the Department of Electrical and Electronic
Engineering, University of Cagliari,  Italy. e-mail: seatzu@diee.unica.it.
\newline
$\-$ $\-$ $\-$   Cort\'es is with the
Department of Mechanical and Aerospace Engineering,
University of California, San Diego.
e-mail:
cortes@ucsd.edu}.
\thanks{
Seatzu's  work is supported by the Region of Sardinia Project RASSR05871 MOSIMA,
 FSC 2014-2020, Annuity 2017, Subject Area 3,
Action Line 3.1.\newline
 $\-$ $\-$ $\-$    Cort\'es' work is supported by NSF Award CNS-1446891.
     \newline
$\-$ $\-$ $\-$  Egerstedt's  work is supported by grant DCIST CRA W911NF-17-2-0181 from the US Army Research Lab.
}}
\maketitle

\begin{abstract}
This paper presents a control technique for output tracking of reference signals in continuous-time
dynamical systems. The technique is comprised of the following three elements: (i) output prediction which has to track the reference
signal, (ii) a controller based on an integrator with variable gain, and (iii) a speedup of the control action
for enhancing the tracker's accuracy and, in some cases, guaranteeing stability of the closed-loop system. The technique is suitable
for  linear and nonlinear systems, implementable by simple algorithms,    can track reference points as well as time-dependent reference signals, and may have large, even global domains of attraction.
 The derived theoretical results
 include  convergence of the tracking controller  and error analysis, and are supported by illustrative
  simulation and laboratory experiments.
\end{abstract}


\section{Introduction}
The subject of this paper is a   reference-tracking control technique for dynamical systems modelled by ordinary differential equations. The technique is founded on real-time implementations of
 a
 fluid-flow variant of the Newton-Raphson method for solving  algebraic equations. The
relevance of the Newton-Raphson method is due to the observation,
argued for in the sequel, that  tracking can be viewed as a dynamic process of attempting to solve a time-dependent suite of nonlinear algebraic equations.

Existing nonlinear regulation techniques   such as the Byrnes-Isidori regulator \cite{Isidori90},
 Khalil's high-gain observers for output regulation \cite{Khalil98},  and
 Model Predictive Control (MPC) \cite{Rawlings17} are  more general and perhaps more powerful than the technique
presented here.       However, their effectiveness is partly due to
 significant computational sophistication like nonlinear
inversions, the appropriate nonlinear normal form, and real-time algorithms for optimal control. The  control technique described in this paper, essentially comprising
a variable-gain integrator,
 is simple and  requires low computing efforts. Nevertheless it will be shown to have inherent stability properties and work
 well on various test problems. Furthermore, it is not a local method, but its domain of attraction is
 often large and sometimes global. Nor is it based on a linearization, and it can be nonlinear. As a matter of fact, the controller is not defined by
 an explicit algebraic   function of the systems' state variable, but rather by a differential equation. As the purpose of this paper is to introduce  a new idea,  we do not make direct comparisons of the proposed technique with existing nonlinear-control methods.
 Instead, we  describe it  in a general setting,
 analyze  its salient features,
 provide results of simulation and laboratory experiments,  and discuss directions for future developments.

 The system-diagram that we consider is depicted in Figure 1, where the reference signal $r(t)$, control input $u(t)$, and system output $y(t)$ are all in $R^m$, with a given $m\in\{1,2,\ldots,\}$.
  The condition that the reference, control, and output have the same dimension is essential for the discussion here, and although ad-hoc ways to circumvent the effects of its absence have begun to emerge \cite{Shivam19a}, we defer their general exposition to a future publication.

   The plant subsystem in Figure 1  is a dynamical system based
  on an ordinary differential equation, whose input, state, and output variables are $u(t)\in R^m$, $x(t)\in R^n$ for some $n\in\{1,2,\ldots\}$, and $y(t)\in R^m$, respectively. The tracking/regulation technique, implemented by the controller subsystem,  is based on the following three elements: output prediction, Newton-Raphson flow, and controller speedup. The predictor computes, at time $t$, an estimate of the future
  output at time $t+T$ for a given $T>0$,
  denoted by $\hat{y}(t+T)$, and the controller is underscored by a process aiming at solving the time-dependent equation $r(t+T)-\hat{y}(t+T)=0$.\footnote{Details of this will be provided in the sequel.} The predictor $\hat{y}(t+T)$  is a function of  $x(t)$ and the input variable $u(t)$, and  therefore  the tracking controller defines
  $\dot{u}(t)$ via a differential equation in terms of $(x(t),u(t))$ as well.  Under the ideal conditions of perfect
 output prediction, this feedback law results in
  perfect asymptotic tracking under general assumptions. In the presence of prediction errors, the asymptotic tracking error will be shown to be equal to the asymptotic prediction error.
 Furthermore, it will be proved that
   an increase in the controller's gain can, in some cases,  stabilize the closed-loop system  and reduce tracking errors that are due to certain disturbances  and computational errors in the loop.
  All of this will be defined and described in detail in later sections.

  \begin{figure}	
	\centering
	\includegraphics[width=0.75\linewidth]{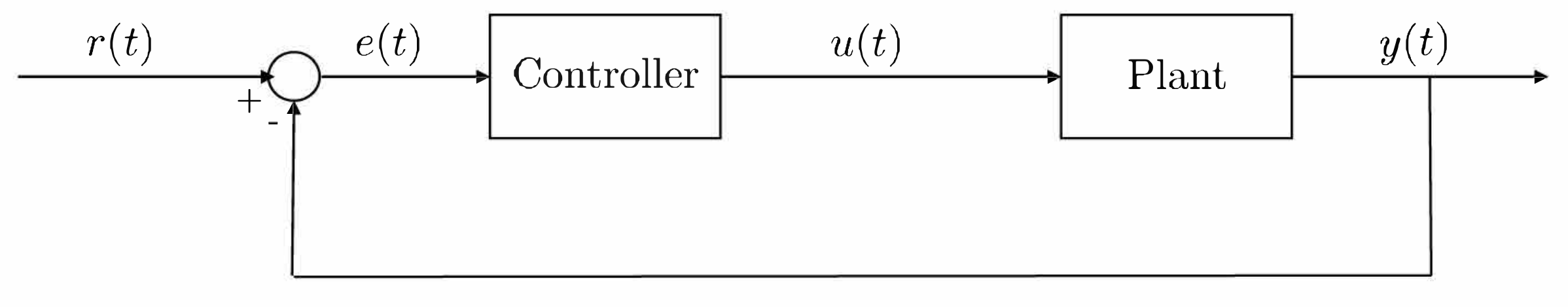}
	\caption{{\small Basic control system}}
\end{figure}

To explain the term ``fluid-flow variant'' of an
iterative algorithm
in $R^m$, and place the forthcoming results in the context of the established literature, consider an iterative algorithm  of the form
\begin{equation}
u_{k+1}=u_{k}+g(u_{k}),
\end{equation}
$k=0,1,\ldots$, where $u_{k}\in R^m$,  and  $g:R^m\rightarrow R^m$ is a function assumed to be locally Lipschitz continuous. Placing the algorithm in a temporal framework, suppose that an iteration according to (1) is computed once every $\Delta t$ seconds for a given $\Delta t>0$, and scale the step size in  the Right-Hand Side (RHS) of (1) by $\Delta t$.    Defining  $u(k\Delta t):=u_{k}$ and
taking the limit  $\Delta t\rightarrow 0$ we obtain the following equation,
\begin{equation}
\dot{u}(t)=g(u(t)).
\end{equation}
The process defined by Eq. (2) is said to be  the {\it fluid-flow version of the algorithm} defined by
Eq. (1).

  Fluid-flow processes can be useful in investigating asymptotic properties of their associated
  discrete algorithms with small step sizes, such as
   convergence,
   optimality and stability of limit points, etc.
   They have been applied mainly to the design of gradient-descent algorithms for problems in optimization and linear algebra, including sorting, eigenvalue decomposition, and linear  programming; see
  \cite{Arrow58,Brockett91,Helmke94}
  for early works. Ref. \cite{Brockett91}
  recognized their potential applications in  massively-parallel computing platforms  such as neural nets,  slated to solve very-large scale problems. Recent applications to learning and  distributed optimization can be found in
  \cite{Lee16,Dhingratoappear}, respectively,  and references therein.

  Second-order optimization methods, especially variants of the Newton-Raphson algorithm, have been considered as well due to their superlinear convergence rates. Refs.
  \cite{Su09,Rahili17} consider first- and second-order algorithms for convex (or concave) constrained programs  with time-varying cost functions. \cite{Su09} is concerned with applications
  to traffic engineering in telecommunications, and \cite{Rahili17}   considers distributed optimization over  multi-agent networks with consensus constraints. Both references derive general theoretical results in abstract  settings of the Newton-Raphson flow beyond their motivating problem-classes, including convergence under weak smoothness assumptions and convergence in a general network setting, respectively.   Ref. \cite{Dhingra17} derives a continuous-flow, primal-dual technique for convex optimization without assuming differentiability of the cost function. Combining results from the theory of convex, nondifferentiable optimization with fluid-flow techniques, it defines the flow by differential inclusions, and derives convergence results, including global asymptotic stability of the minimum and superlinear/quadratic (depending on  assumptions) convergence rates under the weakest-to-date smoothness assumptions on the cost function.

  This paper  also  considers a fluid-flow variant of the Newton-Raphson method, but applies if to finding roots of algebraic equations rather than to convex optimization. However,
  in contrast with the aforementioned references, the resulting control variable $u(t)$ cannot be defined or described by an equation like (2). To get around this difficulty we use an output predictor, $\hat{y}(t+T)$, and define the tracking controller
  as a fluid-flow version of the Newton-Raphson method aiming at solving the time-dependent equations  $r(t+T)-\hat{y}(t+T)=0$.
  Thus, in a way, the predictor is utilized as an observer as well.

  Prediction is commonly used in control, especially for system identification and model construction. In this paper we evidently use it  in a different way, in  the definition of the tight-loop control.
  It resembles the use of prediction in model predictive control
  \cite{Rawlings17}, but our proposed controller is not based on MPC since it does not solve optimal control problems in the loop. A more detailed comparison with MPC will be made at the end of Section 2   following an expanded  explanation of our technique.

  Regarding the third element of the proposed technique, the
   idea that  high controller-gains can enhance stability-robustness  and reduce tracking-convergence time is implicit in    \cite{Sontag89} and explicit in \cite{Kolathaya18}. This paper explores it, in conjunction with the Newton-Raphson flow and output prediction, in the aforementioned general setting of linear and nonlinear control.

The rest of the paper is structured as follows. Section 2 presents the problem and recounts the past developments of our tracking-control technique.
Section
3 carries out analyses of tracking-convergence and how it is impacted by disturbances and other errors in the loop. Section 4  derives a verifiable sufficient condition for stability of linear systems at high controller speeds (rates).
Section 5 presents  simulation results and Section 6 describes a laboratory experiment. Section 7 concludes the paper and
discusses directions for future research.

Preliminary results concerning the material in this paper  can be found in four conference papers, \cite{Wardi17, Wardi18, Shivam19, Shivam19a}.   This paper extends them in the following ways: 1.) It presents a new version of the controller which can yield perfect asymptotic tracking, in contrast with the published works where only approximate tracking is obtained. 2). The analysis includes tracking convergence for general classes of systems and problems (see Section 3). In contrast, the only analyses in the published works pertain to simple examples (see Section 2). 3). It provides a comprehensive stability analysis for  linear systems. 4). The examples are more complex  than in the conference versions.

\section{Problem Formulation and Earlier Results}

This section provides a background material on the specific problem considered in the paper, and recounts the early approaches to it pursued by the authors.

The  output tracking of a system can be viewed as a real-time   implementation of an algorithm for solving a certain kind of algebraic equations.
To see this point consider the system depicted in Figure 1, and
suppose for a moment that the plant-subsystem is a memoryless nonlinearity of the form
\begin{equation}
y(t)=g(u(t)),
\end{equation}
where the function $g:R^m\rightarrow R^m$ is  continuously differentiable.
The tracking problem can be viewed as an effort to solve the time-dependent system of  equations
\begin{equation}
r(t)-g(u)=0
\end{equation}
in the variable $u\in R^m$,
and the controller has to be designed to ensure that
\begin{equation}
\lim_{t\rightarrow\infty}\big(r(t)-g(u(t))\big)=0.
\end{equation}
To solve this problem we define the controller subsystem in Figure 1 so as to implement the fluid-flow version of the Newton-Raphson method. It has the following form,
\begin{equation}
\dot{u}(t)=\Big(\frac{\partial g}{\partial u}(u(t))\Big)^{-1}\big(r(t)-g(u(t))\big),
\end{equation}
where we assume that the Jacobian $\frac{\partial g}{\partial u}(u(t))$ is nonsingular for all $t\geq 0$.

We make the observation that this controller has the form $\dot{u}(t)=A(t)e(t)$, with
$e(t):=r(t)-y(t)$ (see Figure 1) and $A(t)$ is the inverse Jacobian in Eq. (6), an $m\times m$ matrix. This controller  essentially is an integrator with a variable  gain, $A(t)$.
Now it is well known that a standalone integral controller can guarantee tracking of a constant reference, but may result in oscillations and even instability of the closed-loop system (e.g., \cite{Franklin19}). In the present case, however, tracking and stability seem to be ensured  by the particular choice of the gain $A(t)$ in Eq. (6). To see this point, consider the Lyapunov function
\begin{equation}
V(u(t)):=\frac{1}{2}\|r(t)-g(u(t))\|^2.
\end{equation}
In the case where the reference signal is a constant, i.e., $r(t)\equiv r$ for some $r\in R^m$, it can be seen that
$\dot{V}(t)=-\|r-g(u(t))\|^2$, hence by (7)
$\dot{V}(t)=-2V(t)$, which guarantees
asymptotic tracking in the sense of Eq. (5). Variants of
 this simple yet powerful argument underscore  convergence-proofs of
 fluid-flow convex-optimization algorithms in various settings, and the global asymptotic stability of their limit points.  For instance,  see \cite{Dhingra17} for nondifferentiable problems.

If  $r(t)$ is a function of time, Eq. (5) is no longer necessarily true.
However, if $r(t)$ is bounded, continuous, and piecewise
differentiable, the controller defined by
(6) guarantees that
\begin{equation}
{\rm limsup}_{t\rightarrow\infty}\big\|r(t)-y(t)\big\|\leq\eta,
\end{equation}
where $\eta:={\rm limsup}\{\|\dot{r}(t)\|~:~t\in R\}$ (see \cite{Wardi18}).

To tighten  the upper bound in (8), we speed up the action of the controller. One way to do it is to multiply
the Right-Hand Side (RHS) of Eq. (6) by a constant $\alpha>1$, which results in  the following equation,
\begin{equation}
\dot{u}(t)=\alpha\Big(\frac{\partial g}{\partial u}(u(t))\Big)^{-1}\big(r(t)-g(u(t))\big).
\end{equation}
This gives the following bound,
\begin{equation}
{\rm limsup}_{t\rightarrow\infty}\big\|r(t)-y(t)\big\|\leq\frac{\eta}{\alpha},
\end{equation}
provided that the Jacobian $\frac{\partial g}{\partial u}(u(t))$ is nonsingular for every $t\geq 0$
(see \cite{Wardi18}).

This paper considers the plant subsystem to be a dynamical system defined by an ordinary  differential  equation.
Accordingly, let $x(t)\in R^n$ denote its state variable modelled
by  the equation
\begin{equation}
\dot{x}(t)=f(x(t),u(t)),
\end{equation}
where $u(t)\in R^m$ is the control input, $f:R^n\times R^m\rightarrow R^n$ is a suitable function, $t\geq 0$,  and
a given $x(0):=x_{0}\in R^n$ is the initial state. The output function is
\begin{equation}
y(t)=h(x(t)),
\end{equation}
where $y(t)\in R^m$,
for a function $h:R^n\rightarrow R^m$.
We make the following assumptions on the functions $f$ and $h$:
\begin{assumption}\label{ass:ass1}
1). The function $f:R^n\times R^m\rightarrow R^n$ is continuously differentiable, and for every compact set $\Gamma\subset R^m$ there exists $K>0$ such that, for every $x\in R^n$ and for every $u\in\Gamma$,
\begin{equation}
\|f(x,u)\|\leq K\big( \|x\|+1\big).
\end{equation}
2). The function $h:R^n\rightarrow R^m$ is continuously differentiable.
\end{assumption}

Assumption~\ref{ass:ass1} guarantees the existence of a unique continuous, piecewise-differentiable solution for Eq. (11) on the time-horizon  $\{t:t\geq 0\}$,  as long as the input $u(t)$ is
piecewise continuous and bounded.

Extensions of the controller defined in (9) from the case of memoryless plants to that of dynamic plants raises a few challenges. To start with, the input-to-output relation cannot be expressed in a functional algebraic form like in Eq. (3), because $x(t)$, hence $y(t)$ is not a function of $u(t)$ but of $\{u(\tau):\tau<t\}$. Therefore the controller cannot be defined by an equation like (9).
We resolve this issue with the use of an output predictor. Given $T>0$,  it predicts, at time $t$, the future output
 $y(t+T)$, and we denote the predicted value
by $\hat{y}(t+T)$. Suppose that  $\hat{y}(t+T)$ depends on, and is computable from  $x(t)$ and $u(t)$, then it has the following functional form,
\begin{equation}
\hat{y}(t+T):=g(x(t),u(t)),
\end{equation}
where the dependence of $g$ on $T$ is implicit since $T$ is assumed fixed.

\begin{assumption}\label{ass:ass2}
The function $g(\cdot,\cdot)$ is continuously differentiable in $(x,u)$.
\end{assumption}

We define a specific predictor below.  Now the Newton-Raphson flow can extend the one in Eq. (9) by considering  the equation
$r(t+T)-\hat{y}(t+T))=0$ at time $t$. The resulting controller equation has the following
form,
\begin{equation}\label{eq:control_short}
\dot{u}(t)=\alpha\Big(\frac{\partial g}{\partial u}(x(t),u(t))\Big)^{-1}\big(r(t+T)-g(x(t),u(t))\big),
\end{equation}
where it is assumed that $r(t+T)$ is known in advance at time $t$. Putting together the
state equation (11) with the control equation (15), we obtain
the  joint equation
\begin{equation}
\left(
\begin{array}{c}
\dot{x}(t)\\
\dot{u}(t)
\end{array}
\right) ~=
\left(
\begin{array}{c}
 f(x(t),u(t))\\
\alpha\Big(\frac{\partial g}{\partial u}(x(t),u(t))\Big)^{-1}\big(r(t+T)-g(x(t),u(t))\big)
\end{array}
\right),
\end{equation}
which can be viewed as the state equation of an $n+m$-dimensional
dynamical system  with the augmented state $(x(t)^{\top},u(t)^{\top})^{\top}$ and the input $r(t+T)$.
We are concerned with its Bounded-Input-Bounded-State (BIBS)
stability, namely a guarantee that, if the process
 $\big\{r(t+T)~:~t\geq 0\big\}$ is bounded then
 $\{x(t)~:~t\geq 0\}$ and $\{u(t)~:~t\geq 0\}$ are bounded as well. Henceforth we will use the term
 ``stability'' to refer to  BIBS stability.
In contrast with the case where the plant is memoryless, stability cannot be taken for granted in the present case where the plant  is a dynamical system.

This controller  was presented in \cite{Wardi17} with the  particular predictor defined as follows: At time $t$,
given $x(t)$ and $u(t)$, let $\big\{\xi(\tau)~:~\tau\in[t,t+T]\big\}$ be defined by the differential equation
\begin{equation}
\dot{\xi}(\tau)=f(\xi(\tau),u(t)),~~~~~~~\xi(t)=x(t)
\end{equation}
with the boundary condition $\xi(t)=x(t)$; then
 define
\begin{equation}
\hat{y}(t+T):=g(x(t),u(t))=h(\xi(t+T)).
\end{equation}
Observe that Eq. (17) is essentially the state equation (11)
except that it is defined only on the interval $\tau\in[t,t+T]$ with the constant input $u(\tau)\equiv u(t)$ and the initial condition $\xi(t)=x(t)$.

The resulting predictor $\hat{y}(t+T)$, defined by Eqs. (17)-(18),
can admit   efficient  approximations by numerical means such as the Forward Euler method.
Stability of the closed-loop system defined by Eq. (16) with this particular
predictor  was examined (in \cite{Wardi17}) for a number of second-order  linear-system examples. It was shown that,
for a fixed $\alpha$, the system is stable for a large $T$
but unstable for a small $T$. At the same time, small $T$ may be desirable since it
results in a smaller prediction error than larger $T$. To circumvent this conundrum, it was proved that for all of the examples analyzed in \cite{Wardi17}, if the closed-loop system is unstable for  given $T>0$ and $\alpha>0$ then  it can be stabilized by increasing $\alpha$ while keeping the same $T$. Moreover, simulation results suggest that the following  extension of Eq. (10),
\begin{equation}
    {\rm limsup}_{t\rightarrow\infty}\|r(t)-\hat{y}(t)\|<\frac{\eta}{\alpha},
\end{equation}
is satisfied under general conditions.
Thus, a controller's speedup by choosing a large $\alpha$ in Eq. (15)  serves the dual purpose of stabilizing the closed-loop system if need be, and reducing the asymptotic tracking error.
We point out that stabilizability by increasing $\alpha$ is not guaranteed. The derivation of sufficient conditions for it in general
is quite challenging since the function $g(x,u)$ lacks a closed form, but some results will be derived in Section 4.

 Finally, a word must be said about the relationship between the proposed technique and Model-Predictive Control.  MPC uses optimal control over rolling horizons to compute a future target trajectory as well as the  control input  to track it. Our technique is  not concerned with how
 to compute the reference trajectory, but only with its tracking.
  Therefore, if  the reference trajectory is given a priori,  then our technique does not have to solve optimal control problems and it can be simpler than MPC.   On the other hand, if the reference trajectory has to be computed in real time, then our technique can use
  various computational methods  including    interpolation as in Section 5, below; optimal control as
 in MPC or over longer horizons;  or  learning methods based on neural nets as in \cite{Shivam19a}.  Comparisons of its effectiveness and efficiency  vis-a-vis MPC is the subject of a current study.

\section{Enhanced Controller, Tracking and Error Analysis}

This section first presents a modified control algorithm which ensures exact asymptotic tracking of  $r(t+T)$ by $\hat{y}(t+T)$ without resorting to a controller speedup. It then performs an error analysis of the controller which sheds light on the robustness of its tracking performance. In particular, it identifies the errors whose effects on tracking can be reduced by speeding up the controller vs. those whose effects cannot be thus reduced.

In the forthcoming discussion we will use the shorthand notation $\{x(t)\}$ for  the state trajectory $\{x(t):t\geq 0\}$,
and similarly for the input process (trajectory) $\{u(t)\}$, output process $\{y(t)\}$, and
other signals and functions of time.
Also, we will say that the trajectory of the closed-loop system is {\it nonsingular} if for every point
$(x(t)^{\top},u(t)^{\top})^{\top}$ is computes,
the partial Jacobian $\frac{\partial g}{\partial  u}(x(t),u(t))$ is nonsingular.

\subsection{Modified Controller}

Consider the system depicted in Figure 1, where the plant is a dynamical system as defined by Eqs. (11)-(12).  Suppose that Assumption~\ref{ass:ass1} and Assumption 2.2 are satisfied, and $r(\cdot)$ is continuous and piecewise differentiable. Fix a lookahead time $T>0$.  We consider the case where there is no controller speedup, namely $\alpha=1$, and modify the controller equation (15) as follows,

\begin{equation}
\dot{u}(t)=\Big(\frac{\partial g}{\partial u}(x(t),u(t))\Big)^{-1}
\Big(r(t+T) -\hat{y}(t+T)
+\dot{r}(t+T)-\frac{\partial g}{\partial x} (x(t),u(t))f(x(t),u(t))\Big).
\end{equation}

Observe that the difference between this controller and the one defined by (15) is in the addition of the last two terms in the RHS of (20),
$\dot{r}(t+T)-\frac{\partial g}{\partial x} (x(t),u(t))f(x(t),u(t))$.

Define the Lyapunov function
\begin{equation}
    V(x(t),u(t)):=\frac{1}{2}\|r(t+T)-\hat{y}(t+T)\|^2.
\end{equation}

\begin{proposition}\label{prop:prop1}
If the trajectory of the closed-loop system under the state equation (11) and the controller equation (20)  is nonsingular, then $V(x(t),u(t))$ satisfies the following equation,
\begin{equation}
    \dot{V}(x(t),u(t))=-2V(x(t),u(t)).
\end{equation}
Consequently, we have that
\begin{equation}
\lim_{t\rightarrow\infty}\big(r(t)-\hat{y}(t)\big)=0.
\end{equation}
\end{proposition}
{\it Proof.}
Taking the derivative with respect to $t$
in (21), and considering the fact that $\hat{y}(t+T)=g(x(t),u(t))$,
\begin{equation}
\dot{V}(x(t),u(t))=
\big\langle r(t+T)-\hat{y}(t+T),\dot{r}(t+T)-\frac{d}{dt}g(x(t),u(t))\big\rangle.
\end{equation}
Next, by Eqs. (11) and (20),
\begin{eqnarray}
\frac{d}{dt}g(x(t),u(t))
=\frac{\partial g}{\partial x}(x(t),u(t))f(x(t),u(t))
+\frac{\partial g}{\partial u}(x(t),u(t))
\Big(\frac{\partial g}{\partial u}(x(t),u(t))\Big)^{-1}\nonumber \\
\Big(r(t+T)
-g(x(t),u(t))
+\dot{r}(t+T)-\frac{\partial g}{\partial x}(x(t),u(t))f(x(t),u(t))\Big).
\end{eqnarray}
Lastly, simplifying and applying Eq. (25) to (24),  Eq. (22) is obtained.
Consequently, and by (21), Eq. (23) follows.
\hfill$\Box$

\begin{remark}\label{rem:rem1}
Observe that the proof does not require any assumptions about  stability of the closed-loop system. In fact,  if $\{r(t)\}$ is bounded then (by (23)) $\{\hat{y}(t)\}$ is bounded as well, but it is still possible that
$\{||y(t)||\}$ is unbounded. This situation can arise, for example, if the closed-loop system is unstable.
\end{remark}

\subsection{Error Analysis}
This subsection considers three  types of potential errors in the loop, corresponding to the various  terms in the RHS of Eq. (20), and  evaluates their effects on the tracking performance.

\subsubsection{
 Prediction error.} Consider  a prediction error defined as
 $\mathcal{ E}_{1}(t):=\hat{y}(t+T)-y(t+T)$,
 and define the asymptotic prediction error by
\begin{equation}
\eta_{1}:={\rm limsup}_{t\rightarrow\infty}\|{\mathcal E}_{1}(t)\|.
\end{equation}
By Eq. (23),
\begin{equation}
    {\rm limsup}_{t\rightarrow\infty}\|r(t)-y(t)\|=\eta_{1}.
\end{equation}
Defining the asymptotic tracking error by  the Left-hand Side (LHS) of (27),    we see  that
 the asymptotic prediction error is translated to the asymptotic tracking error.

\subsubsection{Error in
$\dot{r}(t+T)-\frac{\partial g}{\partial x}(x(t),u(t))f(x(t),u(t))$.}
Let ${\mathcal E}_{2}(t)$ denote an additive error in the term
$\dot{r}(t+T)-\frac{\partial g}{\partial x}(x(t),u(t))f(x(t),u(t))$
in the RHS of Eq. (20).
Due to this error the controller's definition is modified from (20) to the following equation,
\begin{eqnarray}
\dot{u}(t)=\Big(\frac{\partial g}{\partial u}(x(t),u(t))\Big)^{-1}
\Big(r(t+T) -\hat{y}(t+T)\nonumber \\
+\dot{r}(t+T)-\frac{\partial g}{\partial x} (x(t),u(t))f(x(t),u(t))+{\mathcal E}_{2}(t)\Big).
\end{eqnarray}
Define
\begin{equation}
    \eta_{2}:={\rm limsup}_{t\rightarrow\infty}\|{\mathcal E}_{2}(t)\|.
\end{equation}

\begin{proposition}\label{prop:prop2}
Consider the closed-loop system defined by Eqs. (11), (12), and (28), and suppose that Assumption~\ref{ass:ass1} and Assumption 2.2 are  satisfied. If the trajectory of the closed-loop system   is nonsingular, then
\begin{equation}
    {\rm limsup}_{t\rightarrow\infty}\|r(t)-\hat{y}(t)\|\leq\eta_{2}.
\end{equation}
\end{proposition}

The proof follows as a corollary of Proposition~\ref{prop:prop3} below, hence it is not proved here.

This result, together with the definition of
$\eta_{1}$ (Eq. (26)), imply that
\begin{equation}
    {\rm limsup}_{t\rightarrow\infty}\|r(t)-y(t)\|\leq\eta_{1}+\eta_{2}.
\end{equation}
We next show that  it is possible to reduce the upper bound on the asymptotic tracking error in Eq. (31) by speeding up the controller.  Fix $\alpha>1$, and extend  the definition of  the controller from Eq. (28) to the following equation,
\begin{eqnarray}\label{eq:control_full}
\dot{u}(t)=\Big(\frac{\partial g}{\partial u}(x(t),u(t))\Big)^{-1}
\Big(\big(\alpha\big(r(t+T) -\hat{y}(t+T)\big)\nonumber \\
+\dot{r}(t+T)-\frac{\partial g}{\partial x} (x(t),u(t))f(x(t),u(t))+{\mathcal E}_{2}(t)\Big).
\end{eqnarray}
Observe that the gain $\alpha$ does not multiply the entire RHS of Eq. (32) but only the term
$\big(r(t+T)-\hat{y}(t+T)\big)$ therein.
The result, formalized by the next proposition and the ensuing corollary, shows that it is possible to attenuate the effect of
$\eta_{2}$ but not $\eta_{1}$.

\begin{proposition}\label{prop:prop3}
Consider the closed-loop system defined by Eqs. (11), (12), and (32), and suppose that Assumption~\ref{ass:ass1} and Assumption 2.2 are  satisfied. If the trajectory of the closed-loop system   is nonsingular, then
\begin{equation}
    {\rm limsup}_{t\rightarrow\infty}\|r(t)-\hat{y}(t)\|\leq\frac{\eta_{2}}{\alpha}.
\end{equation}
\end{proposition}
{\it Proof.} Define the Lyapunov function $V(x(t),u(t))$  by Eq. (21). Taking derivatives with respect to $t$, and recalling that
$\hat{y}(t+T)=g(x(t),u(t))$, we have that
\begin{equation}
    \dot{V}(x(t),u(t))=
    \big\langle r(t+T)-\hat{y}(t+T),\dot{r}(t+T)-\frac{d}{dt}g(x(t),u(t))\big\rangle.
\end{equation}
By Eqs. (11) and (32), after some algebra we obtain that
\begin{eqnarray}
\frac{d}{dt}g(x(t),u(t))=\frac{\partial g}{\partial x}(x(t),u(t))f(x(t),u(t))\nonumber \\
+\alpha\big(r(t+T)-\hat{y}(t+T)\big)+\dot{r}(t+T) -\frac{\partial g}{\partial x}(x(t),u(t))f(x(t),u(t))
+{\mathcal E}_{2}(t).
\end{eqnarray}
Using Eq. (35) in Eq. (34) we obtain,
\begin{equation}
    \dot{V}(x(t),u(t))=\big\langle r(t+T)-\hat{y}(t+T),
-\alpha\big(r(t+T)-\hat{y}(t+T)\big)-{\mathcal E}_{2}(t)\big\rangle.
\end{equation}
Consequently, for every $\epsilon>0$ and $t\geq 0$,
 if $\alpha\|r(t+T)-\hat{y}(t+T)\|>\|{\mathcal E}_{2}(t)\|+\epsilon$
then, by the Cauchy-Schwarz inequality,  $\dot{V}(x(t),u(t))<-\epsilon||r(t+T)-\hat{y}(t+T)||$. This, together with the definition of $\eta_{2}$ (Eq. 29)), implies Eq. (33) thereby completing the proof. \hfill$\Box$

\begin{corollary}\label{cor:cor1}
Under the conditions of Proposition~\ref{prop:prop3},
\begin{equation}
    {\rm limsup}_{t\rightarrow\infty}\|r(t)-y(t)\|\leq\eta_{1}+\frac{\eta_{2}}{\alpha}.
\end{equation}
\end{corollary}

{\it Proof.} It follows immediately from Proposition~\ref{prop:prop3} and the definition of $\eta_{1}$. \hfill$\Box$

The enhanced controller, defined by Eq. (32),
seems to have better convergence than the earlier controller defined by Eq. (15). However, the latter controller still has a place since it is simpler, and also can be more practical  in situations where $r(t+T)$ is computed in real time (at time $t$) but $\dot{r}(t+T)$ cannot be computed at that time. An intermediate control algorithm between (15) and (32), defined by Eq. (38), is also possible.
\begin{equation}\label{eq:control_full}
\dot{u}(t)=\Big(\frac{\partial g}{\partial u}(x(t),u(t))\Big)^{-1}
\Big(\alpha\big(r(t+T)
-\hat{y}(t+T)\big)
-\frac{\partial g}{\partial x} (x(t),u(t))f(x(t),u(t))\Big).
\end{equation}
For the purpose of analysis, the controllers based on Eqs. (15) and (38) can be viewed as special cases of the controller defined by (32) by setting
${\mathcal E}_{2}(t):=-\dot{r}(t+T)+\frac{\partial g}{\partial x}(x(t),u(t))$, and  ${\mathcal E}_{2}(t)=-\dot{r}(t+T)$, respectively.

\subsubsection{Error in $\Big(\frac{\partial g}{\partial u}(x(t),u(t))\Big)^{-1}$.}
Convergence of the standard Newton-Raphson method for solving nonlinear equations is known to be robust
to errors in the computation of the inverse-Jacobian (see, e.g.,
\cite{Lancaster66}). A similar robustness holds for convergence of the controller defined by Eq. (32) with respect to errors in the term
$\Big(\frac{\partial g}{\partial u}(x(t),u(t))\Big)^{-1}$,
and Eq. (33) still holds if such errors are small enough.
Therefore we henceforth implicitly assume  that the inverse-Jacobian in Eq. (32) is exact.

\section{Stability Analysis}

The experience with simulation examples in   \cite{Wardi17}  suggests that  an increasing of the controller rate $\alpha$ can stabilize the closed-loop system. This motivates us to explore  verifiable conditions under which this happens.  It may be a difficult problem  for general nonlinear systems, because    the controller $u(t)$ is defined implicitly by a differential equation, whose RHS is not explicit but contains a term, $g(x(t),u(t))$, which also is defined by a differential equation. Therefore, while the problem is posed in a general setting, we  carry out an analysis  only for  linear systems and defer the general case for a future study.

Consider a closed-loop system defined by Eqs. (11)-(12), with the controller defined by either (15), (32) with ${\mathcal E}_{2}(t)\equiv 0$, or (38), with a fixed $\alpha>0$.
It  can be viewed as a dynamical system with state variable $(x(t)^{\top},u(t)^{\top})^{\top}\in R^{n+m}$ and input $r(t)\in R^m$. We  call the state $(x(t)^{\top},u(t)^{\top})^{\top}\in R^{n+m}$ the {\it augmented state}, and denote it by $z(t)$.
The input $\{r(t)\}$ is assumed to be a continuous and piecewise continuously-differentiable function of $t$, and we denote the  $L^{\infty}$ norms of $\{r(t)\}$ and $\{\dot{r}(t)\}$  by $\|r\|_{\infty}$ and $\|\dot{r}\|_{\infty}$, respectively. Assume a given compact set $\Gamma\subset R^{n+m}$ such that the initial (augmented) state $z_{0}:=z(0)$ is constrained to $\Gamma$. The stability notion we have in mind is the following  variant of the  concept of BIBS stability, uniform in $\alpha$:

\begin{definition}\label{def:definition2}
The system is $\alpha$-stable if there exist $\bar{\alpha}\geq 0$ and three class-${\mathcal K}$ functions, $\beta(s)$, $\gamma_{1}(s)$ and $\gamma_{2}(s)$ such that, for every initial state $z_{0}\in\Gamma$, input $\{r(t)\}$,   and $\alpha\geq\bar{\alpha}$,
\begin{equation}
    \|z(t)\|\leq\beta(\|z(0)\|)+\gamma_{1}(\|r\|_{\infty})+\gamma_{2}(\|\dot{r}\|_{\infty}).
\end{equation}
\end{definition}
Note the fact that the three class-${\mathcal K}$ functions are independent of $\alpha\in[\bar{\alpha},\infty)$.

The following result ascertains that $\alpha$-stability
implies asymptotic tracking of $r(t)$ by $\hat{y}(t)$.

\begin{proposition}\label{prop:prop4}
Consider the closed-loop system defined by Eqs. (11)-(12) with the controller defined by either (15), (32) with ${\mathcal E}_{2}(t)\equiv 0$, or (38). Suppose  that Assumption~\ref{ass:ass1} and Assumption 2.2 are  satisfied. If the system is $\alpha$-stable then, for every input $\{r(t)\}$ such that
$\|r\|_{\infty}<\infty$ and $\|\dot{r}\|_{\infty}<\infty$,  for every $z(0)\in\Gamma$, and for every nonsingular trajectory $\{z(t)\}$,
\begin{equation}
    \lim_{\alpha\rightarrow\infty}{\rm limsup}_{t\rightarrow\infty}\|r(t)-\hat{y}(t)\|=0.
\end{equation}
\end{proposition}

{\it Proof.}
Consider first the case where the controller is defined by Eq. (32) with
${\mathcal E}_{2}(t)\equiv 0$. Then for every $\alpha>0$, $\eta_{2}=0$, and hence, by Proposition 3.4,
$\lim_{t\rightarrow\infty}||r(t)-\hat{y}(t)||=0$, this implies (40). Next,
consider  the case where the controller is defined by Eq. (15). It is a special case of Eq. (32) with ${\mathcal E}_{2}(t)=-\dot{r}(t+T)+\frac{\partial g}{\partial x}(x(t),u(t))f(x(t),u(t))$. Therefore, if the system is $\alpha$-stable then there exists $\bar{\eta}_{2}>0$ and $\bar{\alpha}>0$ such that, for every $\alpha\geq\bar{\alpha}$, $\eta_{2}\leq\bar{\eta}_{2}$. Now Eq. (40) follows from Eq. (33). Finally, the case where the controller is defined by Eq. (38) is simpler since it corresponds to (32) with ${\mathcal E}_{2}(t)=-\dot{r}(t+T)$.\hfill $\Box$

Consider now the special case where the system is linear and time invariant. Accordingly, it is    defined  by the  equations
\begin{equation}
\dot{x}(t)=Ax(t)+Bu(t),~~~~~~~y(t)=Cx(t),
\end{equation}
 where $A\in R^{n\times n}$, $B\in R^{n\times m}$, and $C\in R^{m\times n}$
are given matrices. Suppose that the controller is defined by either Eq. (15), (32) with ${\mathcal E}_{2}(t)\equiv 0$, or (38); in either case Assumption 2.1 and Assumption 2.2 are satisfied. The respective analyses of these three cases are almost identical,   hence we perform a detailed  analysis only for the case of (15) and point out in context the required modifications for the two other cases.
Furthermore, to simplify the exposition, we assume that $A$ is nonsingular.

Fix $T>0$.
By Eqs. (17)-(18), we have that
\begin{equation}
    g(x(t),u(t))=Ce^{AT}x(t)+CA^{-1}(e^{AT}-I)Bu(t),
\end{equation}
where $I$ denotes the identity matrix. Therefore,
\begin{equation}
    \frac{\partial g}{\partial x}(x(t),u(t))=Ce^{AT},
    \end{equation}
and
\begin{equation}
 \frac{\partial g}{\partial u}(x(t),u(t))=CA^{-1}(e^{AT}-I)B.
\end{equation}
We assume that  the matrix $\frac{\partial g}{\partial u}(x(t),u(t))=CA^{-1}(e^{AT}-I)B$ is nonsingular.

With the controller defined by (15), the closed-loop system
has the form of Eq. (16).
 By Eqs. (15) and (42)-(44)  the controller has the following
form,
\begin{equation}
\dot{u}(t)=\alpha\Big(CA^{-1}(e^{AT}-I)B\Big)^{-1}r(t+T)
-\alpha\Big(CA^{-1}(e^{AT}-I)B\Big)^{-1}Ce^{AT}x(t)~-\alpha u(t).
\end{equation}
Therefore  Eq. (16)  assumes the form
\begin{equation}
\left(
\begin{array}{c}
\dot{x}(t)\\
\dot{u}(t)
\end{array}
\right)
~=
~\Phi_{\alpha}
\left(
\begin{array}{c}
x(t)\\
u(t)
\end{array}
\right) ~
+~\Psi_{\alpha}r(t+T),
\end{equation}
where $\Phi_{\alpha}$ is an $(n+m)\times(n+m)$ matrix having the following block
structure,
\begin{equation}
\Phi_{\alpha}~=
~\left(
\begin{array}{cc}
A & B\\
-\alpha\Big(CA^{-1}(e^{AT}-I)B)\Big)^{-1}Ce^{AT} & -\alpha I
\end{array}
\right),
\end{equation}
and  $\Psi_{\alpha}$ is an $(n+m)\times n$ matrix of the form
\begin{equation}
\Psi_{\alpha}~=~
\left(
\begin{array}{c}
0\\
\alpha\Big(CA^{-1}(e^{AT}-I)B\Big)^{-1}
\end{array}
\right),
\end{equation}
where the block of zeros is $n\times n$.

Observe that $\alpha$ multiplies the last $m$ rows of $\Phi_{\alpha}$ but none of its first $n$ rows, and hence
 we can write $\Phi_{\alpha}$ in the following way,
\begin{equation}
\Phi_{\alpha}=\left (
\begin{array}{cccccc}
\phi_{1,1} & \phi_{1,2} & \cdot & \cdot & \cdot & \phi_{1,n+m}\\
\phi_{2,1} & \phi_{2,2} & \cdot & \cdot & \cdot & \phi_{2,n+m}\\
\cdot & \cdot & \cdot & \cdot & \cdot & \cdot\\
\cdot & \cdot & \cdot & \cdot & \cdot & \cdot\\
\cdot & \cdot & \cdot & \cdot & \cdot & \cdot\\
\phi_{n,1} & \phi_{n,2} & \cdot & \cdot & \cdot & \phi_{n,n+m}\\
\alpha\phi_{n+1,1} & \alpha\phi_{n+1,2} & \cdot & \cdot & \cdot & \alpha\phi_{n+1,n+m}\\
\cdot & \cdot & \cdot & \cdot & \cdot & \cdot\\
\cdot & \cdot & \cdot & \cdot & \cdot & \cdot\\
\alpha\phi_{n+m,1} & \alpha\phi_{n+m,2} & \cdot & \cdot & \cdot & \alpha\phi_{n+m,n+m}
\end{array}
\right)
\end{equation}
for some scalars $\phi_{j,i}$, $j=1,\ldots,n+m$; $i=1,\ldots,n+m$.
The determinant of $sI-\Phi_{\alpha}$ is a two-dimensional polynomial in $\alpha$ and $s$, which we denote by
$P_{\alpha}(s)$. The  standard formula for computing determinants reveals the following result, whose proof can be found in  the appendix.

\begin{lemma}\label{le:lemma1}
For every $i=1,\ldots,m$ there exists a  polynomial $P_{i}(s)$ in $s$, of degree no more than $n+i$, such that,
\begin{equation}
P_{\alpha}(s)
=\sum_{i=0}^{m}\alpha^i P_{m-i}(s).
\end{equation}
\end{lemma}

\begin{remark}\label{rem:rem4}
For the cases where the controller is defined by either  (32) with ${\mathcal E}_{2}(t)\equiv 0$ or (38),
 the only resulting difference to $\Phi_{\alpha}$ is that the entries of
its last $m$ rows are first-order polynomials in $\alpha$ with possibly-nonzero free coefficients (currently they are first-order polynomials whose free coefficients are 0).
That would not affect the validity of Lemma~\ref{le:lemma1} or the  rest of the analysis in this section.
\end{remark}

Since by assumption $\deg(P_{i})\leq n+i$, we can write $P_{i}(s)$ as
\begin{equation}
P_{i}(s)=\sum_{j=0}^{n+i}a_{i,j}s^j
\end{equation}
 for some coefficients $a_{i,j}$, $j=0,\ldots,n+i$.
 We assume, without loss of generality,  that  $a_{i,n+i}\neq 0$ to ensure that ${\rm deg}(P_{i})=n+i$.
 Then
 \begin{equation}
 P_{\alpha}(s)=\sum_{i=0}^{m}\alpha^i\sum_{j=0}^{n+m-i}a_{m-i,j}s^j.
 \end{equation}
 The highest-order term (in $s$) of $P_{\alpha}(s)$ is $a_{m,n+m}s^{n+m}$, and we assume that
 $a_{m,n+m}=1$.

We next derive a sufficient condition for the $\alpha$-stability  of the system.  The condition
consists of two polynomials having all of their roots in the Left-Half Plane (LHP).   One polynomial has degree $n$, the other has degree $m$,
 and both are independent of $\alpha$ hence the sufficient condition is verifiable.

 The first polynomial is $P_{0}(s)$, which by (50) is the polynomial-coefficient of $\alpha^m$, the leading term in $P_{\alpha}(s)$ in terms
 of the power of $\alpha$. Note (Eq. (51)) that $\deg(P_{0})=n$.

 The second polynomial, denoted by $Q(s)$,  is defined as follows. For every $i=0,\ldots,m$,
 consider the polynomial $P_{i}(s)$, defined in Eq. (51), whose degree is $n+i$. Define a polynomial $\tilde{P}_{i}(s)$ as the monomial consisting of the highest-order term of $P_{i}(s)$, namely,
 \begin{equation}
 \tilde{P}_{i}(s)=a_{i,n+i}s^{n+i}.
 \end{equation}
 Next, in analogy to (50),     define the
 family of polynomials parameterized by $\alpha>0$,
 $\{\tilde{P}_{\alpha}(s)\}$, by
 \begin{equation}
 \tilde{P}_{\alpha}(s)=\sum_{i=0}^{m}\alpha^i\tilde{P}_{m-i}(s).
 \end{equation}
 By (53),
 \begin{equation}
 \tilde{P}_{\alpha}(s)=\sum_{i=0}^{m}\alpha^i a_{m-i,n+m-i}s^{n+m-i}.
 \end{equation}
 Observe that for every $\alpha>0$, $\tilde{P}_{\alpha}(s)$ is evenly divisible by $s^n$. Dividing it by $s^n$, we define
 \begin{equation}
 \tilde{Q}_{\alpha}(s):=\sum_{i=0}^{m}\alpha^i a_{m-i,n+m-i}s^{m-i},
 \end{equation}
 and we note that
 \begin{equation}
 \tilde{P}_{\alpha}(s)=s^{n}\tilde{Q}_{\alpha}(s).
 \end{equation}
We make the observation  that $\tilde{P}_{\alpha}(s)$ has the degree (in $s$) of $n+m$ hence it has $n+m$
  roots; by (57), $n$ of those roots are at $s=0$, and the remaining $m$ roots are the roots of $\tilde{Q}_{\alpha}(s)$. Finally, we define the $mth$-degree
  polynomial $Q(s)$ by setting $\alpha=1$ in
  $\tilde{Q}_{\alpha}(s)$ (Eq. (56)); namely,
 \begin{equation}
 Q(s):=\tilde{Q}_{1}(s)=\sum_{i=0}^{m}a_{m-i,n+m-i}s^{m-i}.
 \end{equation}
 Observe that $Q(s)$ is independent of $\alpha$, and its degree is $m$.

 The following result establishes the $\alpha$-stability of the system.

 \begin{theorem}\label{th:theorem1}
 If the polynomials $P_{0}(s)$ and $Q(s)$ have all of their roots in the open Left-Half Plane (LHP), then the system is $\alpha$-stable.
 \end{theorem}


The  proof is based on the following two arguments: For large-enough $\alpha$, (i)    the matrix $\Phi_{\alpha}$ is Hurwitz, and (ii) the effect of the gain   $\alpha$ in $\Psi_{\alpha}$ (Eq. (48)) on  $||z(t)||$ is bounded even though $\alpha$ can be arbitrarily large.

To prove the first argument  we employ a root-locus technique in a nonstandard setting, where  the functional  dependence of $P_{\alpha}(s)$ on $\alpha$ and $s$ is via a two-dimensional polynomial. The proof proceeds as follows:
First we show that bounded branches of the root locus must converge to the zeros of $P_{0}(s)$, and
 this follows standard root-locus arguments. Then we prove that unbounded branches have the same asymptotic angles  as the angles of the roots of $Q(s)$, hence unbounded branches will be in the LHP for large-enough $\alpha$ if all of the roots of $Q(s)$ are in the LHP.

 The  proof of Theorem~\ref{th:theorem1} will be preceded by a sequence of technical lemmas, where those proofs that are straightforward  are relegated to the appendix.  Throughout the forthcoming discussion we denote a generic branch of the root locus of $P_{\alpha}(s)$ by
 $\{s(\alpha)\}_{\alpha\geq 0}$, or  by $\{s(\alpha)\}$ for a
 simpler notation.

\begin{lemma}\label{le:lemma2}
If $\{s(\alpha)\}$ is bounded over $\alpha\in[0,\infty)$,
then the limit
 $
 \lim_{\alpha\rightarrow\infty}s(\alpha)$ exists and it is a root of $P_{0}(s)$.
 \end{lemma}
 For a proof, please  see the appendix.

 Consider next the case where $\{s(\alpha)\}$ is unbounded. Let $A\subset [0,\infty)$ be an unbounded  set such that
 \[
 \lim_{\alpha\rightarrow\infty;~\alpha\in A}|s(\alpha)|=\infty.
 \]

\begin{lemma}\label{le:lemma3}
There exist constants $c>0$ and $C>c$ such that, as $\alpha\rightarrow \infty;~\alpha\in A$,
 \begin{equation}
 c\leq ~{\rm liminf}~\frac{|s(\alpha)|}{\alpha},~~~~~{\rm and}~~~~~{\rm limsup}~\frac{|s(\alpha)|}{\alpha}~\leq C.
 \end{equation}
 \end{lemma}
 {\it Proof.}
Consider first the right inequality of Eq. (59). Let us argue by contradiction. If that inequality does not
 hold,  there exists an unbounded set   $A_{1}\subset A$ such that, as $\alpha\rightarrow\infty$, $\alpha\in A_{1}$,
 \begin{equation}
 \frac{|s(\alpha)|}{\alpha}\rightarrow\infty.
 \end{equation}
 By Eq. (50), $\forall\alpha\in A_1$,
 \[
 \sum_{i=0}^{m}\alpha^i P_{m-i}(s(\alpha))=0.
\]
 Dividing this equation by $s(\alpha)^{m+n}$, we get that
 \begin{equation}
 \sum_{i=0}^{m}\Big(\frac{\alpha}{s(\alpha)}\Big)^i \times\frac{P_{m-i}(s(\alpha))}{s(\alpha)^{n+m-i}}
 =~\sum_{i=1}^{m}\Big(\frac{\alpha}{s(\alpha)}\Big)^i \times \frac{P_{m-i}(s(\alpha))}{s(\alpha)^{n+m-i}}+\frac{P_{m}(s(\alpha))}{s(\alpha)^{m+n}}~=~0.
 \end{equation}
 But deg$(P_{m-i})=n+m-i$, hence, and by (51), as $\alpha\rightarrow\infty;~\alpha\in A_1$,
  \[
  \frac{P_{m-i}(s(\alpha))}{s(\alpha)^{n+m-i}}\rightarrow a_{m-i,n+m-i}
  \]
  which is a finite-magnitude number. Therefore, and by (60),
 \[
 \sum_{i=1}^{m}\Big(\frac{\alpha}{s(\alpha)}\Big)^i\times  \frac{P_{m-i}(s(\alpha))}{s(\alpha)^{n+m-i}}\rightarrow 0
 \]
 as $\alpha\rightarrow 0;~\alpha\in A_1$.
 Furthermore, deg$(P_{m})=n+m$, hence, and since (by assumption) the leading coefficient of $P_{m}$ is 1,
 \[
 \frac{P_{m}(s(\alpha))}{s(\alpha)^{m+n}}\rightarrow 1.
 \]
 This contradicts (61) thereby ascertaining the right inequality of (59).

 The left inequality of (59) is provable by similar arguments, hence it is relegated to the appendix.      \hfill$\Box$

Given  a  polynomial $P_{\alpha}(s)$ (as defined by (50)) and $\alpha\geq 0$,  we next
 examine the derivatives of $s(\alpha)$
 with respect to the coefficients of $P_{m-i}(s)$, for $i=0,\ldots,m$, as defined by  (51).
 For this purpose we consider all but the leading coefficients, namely $a_{m-i,j}$, $j=0,\ldots,n+m-i-1$. We
 denote these derivatives by $\frac{\partial s(\alpha)}{\partial a_{m-i,j}}$.
 For apparent reasons of notation, we will use $\ell$ and $\nu$ instead of $i$ and $j$ in the following discussion

 \begin{lemma}\label{le:lemma4}
 There exist $r\geq 0$ and $L>0$ such that, if $|s(\alpha)|\geq r$, then for every $\ell=0,\ldots,m$, and for every
 $\nu=0,\ldots,n+m-\ell-1$,
 \begin{equation}
\Big|\frac{\partial s(\alpha)}{\partial a_{m-\ell,\nu}}\Big|\leq L.
 \end{equation}
 \end{lemma}
The proof is carried out in the appendix  by realizing that $P_{\alpha}(s(\alpha))=0$, and taking derivatives of this equation  with respect to $a_{m-\ell,\nu}$.
  We remark that the assertion of Lemma~\ref{le:lemma4}  may not hold true for the case where $\nu=n+m-\ell$, namely for
the leading coefficient of $P_{m-\ell}(s)$.

Recall the definition of $\tilde{P}_{\alpha}(s)$ which was made in Eq. (55). Similarly to the notation $s(\alpha)$ for a generic root of $P_{\alpha}(\cdot)$, we denote
by $\{\tilde{s}(\alpha)\}$ a generic branch of the root locus  of $\tilde{P}_{\alpha}(\cdot)$.

\begin{lemma}\label{le:lemma5}
There exist constants $r>0$ and $K>0$     such that,
if $|s(\alpha)|\geq r$ for some $\alpha>0$, then there exists $\tilde{s}(\alpha)$ such that
\begin{equation}
|\tilde{s}(\alpha)-s(\alpha)|<K.
\end{equation}
\end{lemma}

{\it Proof.} The polynomials $P_{m-\ell}(s)$ and $\tilde{P}_{m-\ell}(s)$,
$\ell=0,\ldots,m$, have the same  respective leading coefficients, $a_{m-\ell,n+m-\ell}$. As for the other coefficients, those of $P_{m-\ell}(s)$
are $a_{m-\ell,\nu}$,
$\nu=0,\ldots,n+m-\ell-1$,
and those of $\tilde{P}_{m-\ell}(s)$ are 0.
The statement now follows
from  Lemma~\ref{le:lemma4} and the mean-value theorem. \hfill $\Box$

Fix $\alpha>0$. It has been mentioned that, by Eq. (57), $n$ of the roots of $\tilde{P}_{\alpha}(s)$ are at 0, and its  remaining $m$  roots are the roots of
$\tilde{Q}_{\alpha}(s)$ as defined by (56). We next characterize the roots of
$\tilde{Q}_{\alpha}(s)$.

\begin{lemma}\label{le:lemma6}
Let $s$ be a root of the polynomial $Q(\cdot)$. Then
for every $\alpha>0$, $\alpha s$ is a root of the polynomial $\tilde{Q}_{\alpha}(\cdot)$.
\end{lemma}

{\it Proof.} By Eqs. (56) and (58), we see that for every complex variable $s$, and for every $\alpha>0$,
\begin{equation}
\tilde{Q}_{\alpha}(\alpha s)=\alpha^m Q(s).
\end{equation}
Therefore, if $s$ is a root of $Q(\cdot)$, $\alpha s$ is a root of $\tilde{Q}_{\alpha}(\cdot)$. \hfill $\Box$

Given a complex variable $s$, let $\angle s$ denote the angle (argument)
of $s$ with respect to the positive side of the horizontal axis. Thus, if
$s=|s|e^{j\phi}$ according to its polar coordinates, then $\angle s=\phi$.

\begin{lemma}\label{le:lemma7}
 Let $s_{i}$, $i=1,\ldots,m$ denote the roots of the polynomial $Q(s)$. Suppose that none of these roots is 0. For every unbounded branch of the root locus of $P_{\alpha}(s)$, denoted by
$\{s(\alpha)\}$, there exists $i\in\{1,\ldots,m\}$ such that,
\begin{equation}
\lim_{\alpha\rightarrow\infty}\angle s(\alpha)=\angle s_{i}.
\end{equation}
\end{lemma}

{\it Proof.} By Lemma~\ref{le:lemma6}, $m$ of the root-locus' branches  of $\tilde{Q}_{\alpha}(s)$
are straight lines $\{\alpha s_{i}\}_{\alpha=0}^{\infty}$,
$i=1,\ldots,m$. By Eq. (57), these are the unbounded root loci of $\tilde{P}_{\alpha}(s)$. Therefore, and by
Lemma~\ref{le:lemma5}, if $\{s(\alpha)\}$ is unbounded, there exist
$r>0$,   $K>0$ and $i\in\{1,\ldots,m\}$ such
that, if $|s(\alpha)|>r$,
then $|s(\alpha)-\alpha s_{i}|<K$. This implies Eq. (65) and completes the proof.\hfill $\Box$

{\it Proof of Theorem~\ref{th:theorem1}.}
Suppose that all of the roots of the polynomials $P_{0}(s)$ and $Q(s)$ are in the LHP. Then Lemma~\ref{le:lemma2} and Lemma~\ref{le:lemma7} imply that there exists
$\alpha\geq 0$
such that $\forall\alpha\geq\bar{\alpha}$,
the closed-loop system matrix $\Phi_{\alpha}$ is Hurwitz.  According the Definition~\ref{def:definition2}, we have to show that the class-${\mathcal K}$ functions $\beta$ and  $\gamma_{1}$ satisfy Eq. (39)
for all large-enough $\alpha$ ($\gamma_{2}(\cdot)$ is irrelevant because
$\{\dot{r}(t)\}$ is not an explicit part of the input).
This is not apparent in light of the fact that the matrix $\Psi_{\alpha}$ has a multiplicative $\alpha$-term (see Eq. (48)). Nonetheless this is true because of the block of zeros in $\Psi_{\alpha}$. We next show this point.

As a matter of notation, we say that a matrix is $O(\alpha^{k})$ for an integer $k$ (possibly nonpositive)
if the highest power of $\alpha$ among all of its elements  is $\alpha^k$. Recall Eq. (49), and note,  that the first $n$ rows of $\Phi_{\alpha}$ do not contain $\alpha$, and the last $m$ rows contain $\alpha$ as a multiplicative factor. Therefore, by Cramer's rule, the first $n$ columns of $(sI-\Phi_{\alpha})^{-1}$ are $O(\alpha^0)$,  and
its  last $m$ columns are $O(\alpha^{-1})$. Denote by $\Phi_{1,\alpha}(s)$ and $\Phi_{2,\alpha}(s)$ the matrices comprised of the first $n$ columns and last $m$ columns of $(sI-\Phi_{\alpha})^{-1}$, respectively. Then  $\Phi_{1,\alpha}(s)$ is $O(\alpha^{0})$, and
$\Phi_{2,\alpha}(s)$ is $O(\alpha^{-1})$. As for $\Psi_{\alpha}$, denote the matrix comprised of its last $m$ rows by $\Psi_{2,\alpha}$. Then (by (48)), $\Psi_{2,\alpha}$ is $O(\alpha^1)$. Now the $r$-to-$z$
(input-to-state) transfer function is
\begin{equation}
 (sI-\Phi_{\alpha})^{-1}\Psi_{\alpha}=
 \left(\begin{array}{cc}\Phi_{1,\alpha}(s) & \Phi_{2,\alpha}(s)
 \end{array}
 \right)
 \left(\begin{array}{c}
 0\\
 \Psi_{2,\alpha}\end{array}\right)=\Phi_{2,\alpha}(s)\Psi_{2,\alpha}.
\end{equation}
Since $\Phi_{2,\alpha}(s)$ is $O(\alpha^{-1})$ and $\Psi_{2,\alpha}$ is $O(\alpha^1)$, $(sI-\Phi_{\alpha})^{-1}\Psi_{2,\alpha}$ is $O(\alpha^0)$. Therefore, and since $\Phi_{\alpha}$ is Hurwitz, there exist $\sigma>0$ and $\bar{\alpha}\geq 0$ such that, for every $\alpha\geq\bar{\alpha}$, the real parts all the poles of the $r$-to-$z$ transfer function
are smaller than  $-\sigma$. This implies the $\alpha$-stability of the closed-loop system.\hfill$\Box$

We remark that if either matrix  $P_{0}(s)$ or $Q(s)$ has a root in the RHP then  the closed-loop system is not $\alpha$-stable.

{\it Example.} The following example is of an   $\alpha$-stable system where the plant subsystem is   neither stable not of a minimum phase. Let
\[
A=\left(\begin{array}{cc}
2 & 1\\
-1 & -1
\end{array}
\right),~~~~~~~B=\left(\begin{array}{c}
0\\
1
\end{array}
\right),~~~~~~~C=\left(\begin{array}{cc}
-10 & 1
\end{array}\right),
\]
and $T=0.25$s. The plant transfer function is
\[
G(s)=\frac{s-12}{s^2-s-1},
\]
which is unstable and not of a minimum phase. Next,
$P_{\alpha}(s)=\big(s^3-s^2-s\big)+\alpha\big(s^2+16.19s+97.18\big)$. Therefore   $P_{0}=s^2+16.19s+97.18$ and
$P_{1}(s)=s^3-s^2-s$, implying that   $Q(s)=s+1$. Both $Q(s)$ and $P_{0}(s)$ have all of their roots in the LHP, hence the system is $\alpha$-stable.

\section{Simulation Experiments}

This section presents simulation results for  two problems, namely an inverted pendulum and a platoon  of autonomous vehicles. For the inverted pendulum we  use the controller defined by Eq. (32). As for the platoon system, we assume that the vehicles' controllers have no a-priori knowledge of $\{\dot{r}(t)\}$ for their respective reference trajectories, therefore we use the controllers defined by Eqs. (15) and (38). We then present the better results of the two, which are with (15).

\subsection{Inverted pendulum}
The considered pendulum  is mounted on a  cart which can move in the two directions of a given line, parameterized by $z\in R$.
 Let $\theta$ denote the angle of the pendulum from its pivot on the cart to the left of the upward-vertical direction.
Thus, if the pendulum is pointed upwards then $\theta=0$, and if it points sideways along the $z$ axis in the positive direction then $\theta=-\pi/2$ rads. Let $M$ and $m$ denote the masses of the cart and pendulum, respectively, and let $\ell$ be the distance from the cart to the pendulum's center of mass. Furthermore, let  $F$ be the force applied to the cart in the positive direction of the $z$ axis, and let $\theta$ be the system's output to be controlled.

This system  generally is four-dimensional with the
state variable $x:=(z,\dot{z},\theta,\dot{\theta})^{\top}$.
However, a  simpler, second-order representation of the  pendulum's motion  can be obtained by making
the following two assumptions: 1). The pendulum consists of a weightless rod and a point mass at its end. 2). There is no friction in the movement of either  cart or  pendulum. In this case, the dynamic equation of the pendulum's motion becomes
\begin{equation}
(M\ell+m\ell\sin^2\theta)\ddot{\theta}
+m\ell\dot{\theta}^2(\sin\theta)(\cos\theta)
+~(M+m)g\sin\theta~=~F\cos\theta;
\end{equation}
see \cite{Wikipedia10}.
This equation provides a state-space representation of the system  where the state variable is $x=(\theta,\dot{\theta})^{\top}$, the input is $u=F$, and the output is $y=\theta$. We chose the following parameters for the simulation:
 $M=1$kg, $m=0.2$kg, $\ell=2$m, and $g=9.81m/s^2$. The simulation starts at the initial
state $x(0)=(\frac{\pi}{6},0)^{\top}$, and it
solves the state equation in a specified horizon
$t\in[0,t_{f}]$ by the forward Euler method with
the step-size $dt=0.01$s. The control algorithm uses the prediction horizon $T=0.2$s, and it computes the predicted state trajectory (Eq. (17)) by the forward-Euler method with the step-size $\Delta t=0.01T$. The initial condition for the controller
equation (32)  is $u(0)=0$.

The target trajectory for the tracking-control experiment is
$r(t)=-\frac{\pi}{6}+\frac{\pi}{3}\sin t$, which oscillates between the angles of $30^o$ and $-90^o$. At $-90^o$ the pendulum points at the horizontal direction along the positive z-axis,  and this can be problematic because it is physically impossible to balance the pendulum at this angle.  However,  in the present experiment $r(t)$ just touches the horizontal direction and then immediately retreats therefrom.
The
time-horizon for the simulation is $t_{f}=25$s.

For the controller's equation (32) we first  took $\alpha=1$, and noted  convergence of $\theta(t)$ to $r(t)$ in about 2 seconds. To speed up the convergence we increased the controller's gain to $\alpha=35$, and the results are depicted in Figures 2-4.
 Figure 2 shows the graphs of $\theta(t)$ in blue, and the reference $r(t)$
in red. The two graphs appear to coalesce for the first time at about $t=1$s, and remain close to each other except for  slight differences when $r(t)\sim-\frac{\pi}{2}$ rads (about -1.57 in the graph). This is not surprising because at such points the pendulum is horizontal. The maximum error, $|r(t)-\theta(t)|$,
for $t\geq 1$ was measured from the graphs at 0.022 radians, or 1.2605 degrees.
To further highlight the discrepancies between $\theta(t)$ and $r(t)$ we plot the angular velocity, $\dot{\theta}(t):=x_{2}(t)$. The results are depicted in Figure 3, and they clearly show a distortion from the sinusoidal form of $\dot{r}(t)$ at points where $\theta\sim -\frac{\pi}{2}$rads. Furthermore, we plot the graph of
the control signal $u(t)$  in Figure 4, and we notice large peaks at the point where $r(t)=-\frac{\pi}{2}$. All of this  is expected in light of the earlier remarks concerning the challenges of controlling the pendulum  at (or close to)  the horizontal angle.

To verify that the discrepancies between $r(t)$ and $\theta(t)$, and the large peaks in $u(t)$  indeed are due to the fact that $r(t)$ reaches $-\pi/2$ periodically,  we attenuated the sinusoid
part of $r(t)$ by the factor of $0.8$, and thus  $r(t)=-\frac{\pi}{6}+0.8\frac{\pi}{3}\sin t$, corresponding to  oscillations between the angles of $18^o$ and $-78^o$. We only show the resulting graph of $\dot{\theta}(t)$ since it most clearly indicates the distortions in $\theta(t)$. This graph is depicted in  Figure 5, where its distortion at about the lower-peak angle of    $-78^o$ is barely visible. The discrepancies between $\theta(t)$ and $r(t)$ is hardly noticeable  from their respective graphs which are not shown here. Also, the peak control at these values (not shown here) is reduced to nearly 60, which is about 10\% of its value obtained from the full-sinusoidal swing that is depicted in Figure 4.

 \begin{figure}	
	\centering
	\includegraphics[width=0.75\linewidth]{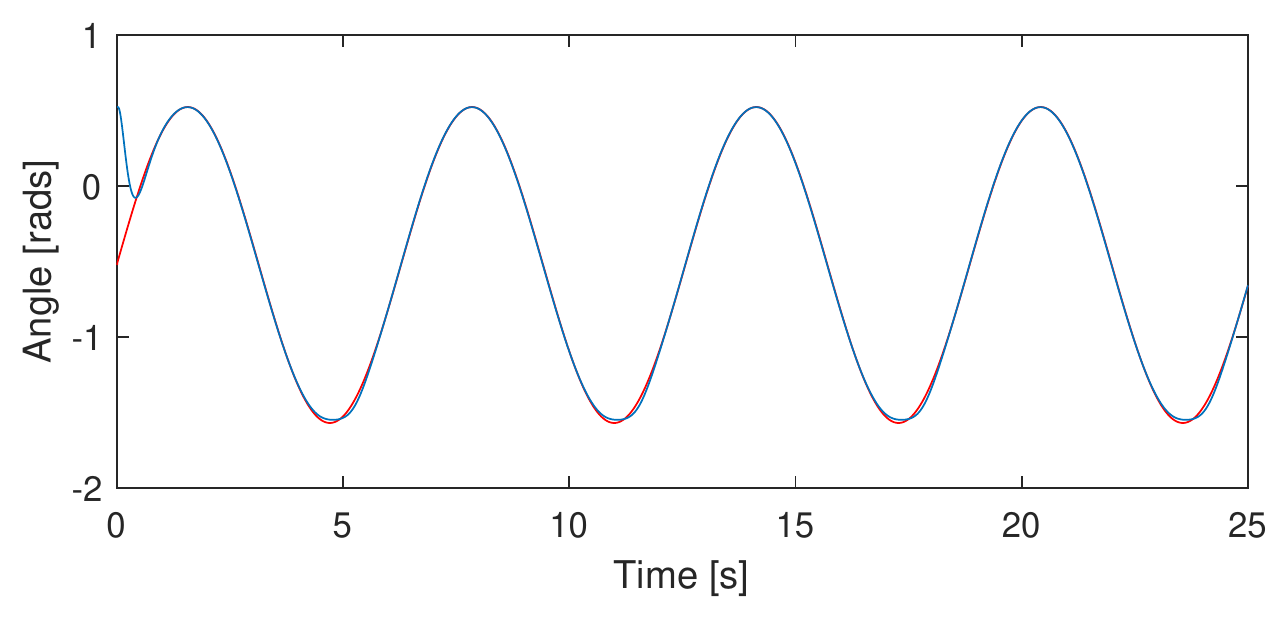}
	\caption{{\small Inverted pendulum: $\theta$  and $r$ vs. $t$}}
\end{figure}

 \begin{figure}	
	\centering
	\includegraphics[width=0.75\linewidth]{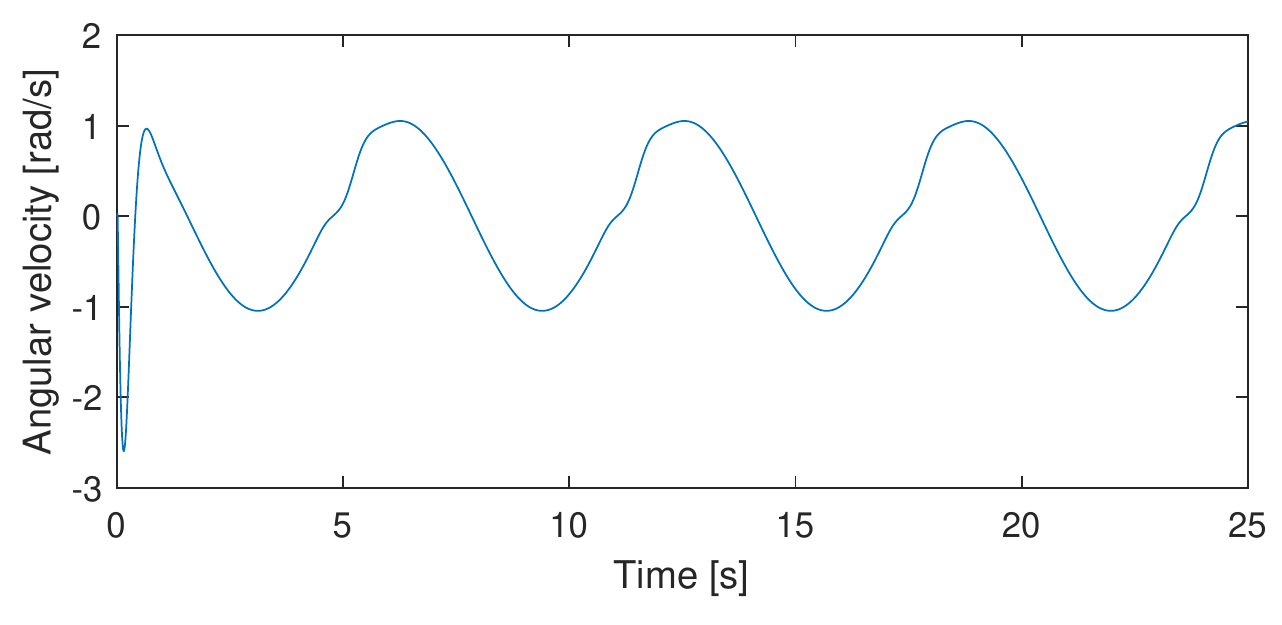}
	\caption{{\small Inverted pendulum: $\dot{\theta}$  vs. $t$}}
\end{figure}

 \begin{figure}	
	\centering
	\includegraphics[width=0.75\linewidth]{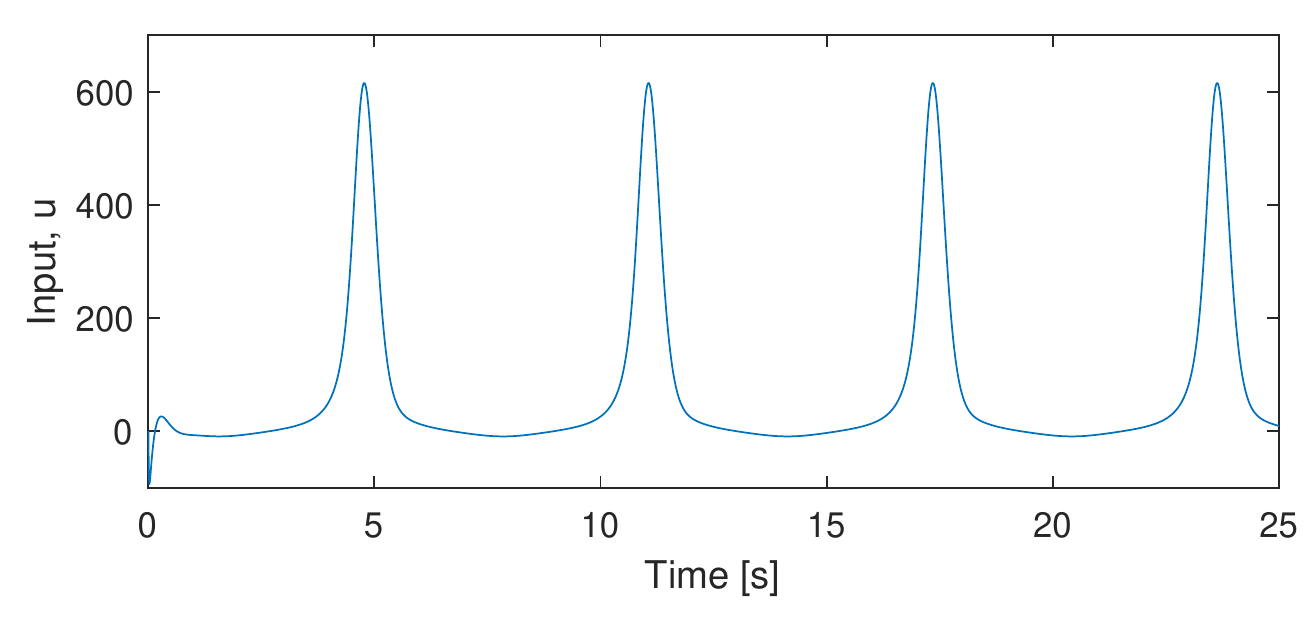}
	\caption{{\small Inverted pendulum: $u$ vs. $t$}}
\end{figure}

 \begin{figure}	
	\centering
	\includegraphics[width=0.75\linewidth]{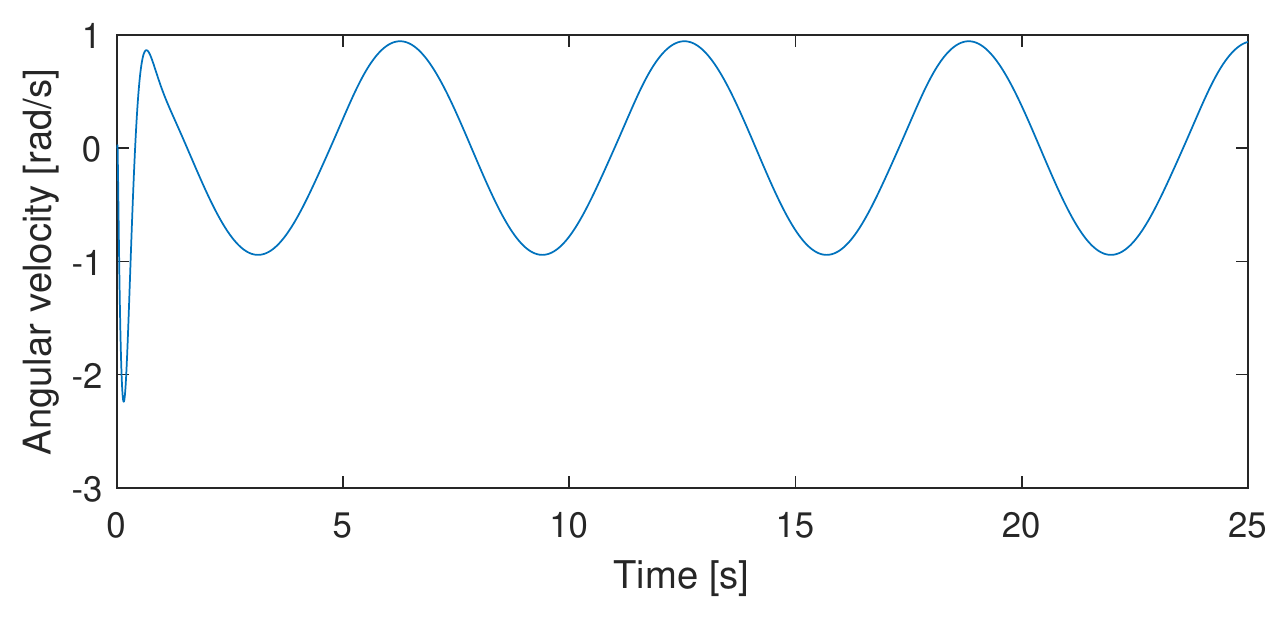}
	\caption{{\small Inverted pendulum, reduced oscillations: $\dot{\theta}$ vs. $t$ }}
\end{figure}

\subsection{Platoon of autonomous vehicles}
The simulation experiment described in this subsection concerns the planar  motion of a platoon, controlled to
follow a given path
 in the $(z_{1},z_{2})$ plane. The platoon consists of four agents (vehicles), denoted by  $A_{i}$, $i=1,2,3,4$, in the order of their movement.  $A_{1}$ is the leading vehicle, and $A_{i}$ follows $A_{i-1}$, $i=2,3,4$.
$A_{1}$ is provided with an exogenous  reference trajectory (path) to track, $\{r_{1}(t)\}$, and for $i=2,3,4$, $A_{i}$ attempts to follow $A_{i-1}$ at a prescribed distance (arclength) of  $d$ m
on the path. Whereas  the target reference for each agent remains on the path $\{r_{1}(t)\}$, the agent itself can get off the path while pursuing its target reference. In this way the agents' motions are two-dimensional and not confined to one-dimensional curves.  We assume that
each agent controls its own motion: $u_{1}(t)$ depends on $r_{1}(t)$, while for   $i=2,3,4$,
 $u_{i}(t)$ is computed by the position and  velocity of $A_{i-1}$, which are assumed to be measured by $A_{i}$ or transmitted to it by $A_{i-1}$.

The motion-dynamics of the vehicles follow the bicycle model, a sixth-order nonlinear system that has been extensively used in the design and analysis of  motion control  for autonomous vehicles;  see, e.g.,  \cite{Plessen18} and references therein.  The state space consists of the  six-tuple
${x=(z_{1},z_{2},v_{\ell},v_{n},\psi,\dot{\psi})^{\top}}$,
where $z_{1}$ and $z_{2}$ are the planer position-coordinates of the center of gravity of the vehicle,  $v_{\ell}$ and  $v_{n}$ are  the longitudinal and lateral velocities,  $\psi$ is the heading of the vehicle and  $\dot{\psi}$ is its angular velocity. The input,  ${u=(a_\ell , \delta_f)^\top}$, consists of  the longitudinal acceleration  and steering angle of the front wheel, respectively, and the output is the position of the center of gravity of the vehicle, namely $y=(z_{1},z_{2})^{\top}$.

The dynamic equations of the vehicles are given by the following equations (see \cite{Kong15}),
\begin{align}
\dot{z}_{1}&=v_{\ell}\cos\psi-v_{n}\sin\psi\\
\dot{z}_{2}&=v_{\ell}\sin\psi+v_{n}\cos\psi\\
\dot{v}_{\ell}&=\dot{\psi}v_{n}+a_{\ell}\\
\dot{v}_{n}&=-\dot{\psi}v_{\ell}+{2}\left(F_{c,f}\cos\delta_f+F_{c,r}\right)/{m}\\
\ddot{\psi}&={2}\left(l_f F_{c,f}\cos\delta_f-l_r F_{c,r}\right)/{I_z},
\end{align}
where $m$ is the mass of the vehicle, $l_f$ and $l_r$ are the  front and rear axles' distances from the vehicle's center of mass,
$I_z$ is the yaw moment of inertia,  and $F_{c,f}$ and $F_{c,r}$ are the lateral forces on the front and rear wheels.  These  forces   are  approximated by the following equations,
\begin{align}
F_{c,f}=C_{\alpha,f}\left(\delta_f-\tan^{-1}\left((v_{n}+l_f\dot{\psi})/v_{\ell}\right)\right)\\
F_{c,r}=-C_{\alpha,r}\tan^{-1}\left(({v_{n}-l_r\dot{\psi}})/{v_{\ell}}\right),
\end{align}
where  $C_{\alpha,f}$ and $C_{\alpha,r}$ are   the cornering stiffness of the front and rear tires, respectively.

In the simulation we used the following model-parameters as in \cite{Mathworks19}, Volvo V70 model, except for $I_{z}$ (not provided there) which has been estimated by averaging data from cars of similar weights and dimensions. Thus,
$m=1,700$kg,
$l_r=1.5$m, $l_f=1.5$m,   $I_z=2,500{\rm kg\cdot m^2}$, and $C_{\alpha_f}=C_{\alpha_r}=29,963.5$N/rad. As for the considered problem, controller and simulation parameters, the desired inter-agent distance is $d=10$m,   the simulation horizon is $t_{f}=38$s,  and the discretization step size for the simulation is $dt=0.01$ secs. The controller prediction horizon is set to $T=0.5$s, and the    discretization time step for the predictor is $\Delta T=0.001$T. The controllers' speedup factor is  $\alpha=100$ for all the vehicles.
The target trajectory $\{r_{1}(t)\}$ is  indicated by the curve in Figure 6, and  its  acceleration along the path is indicated by the blue graph in Figure 7. Its initial speed is $\dot{r}_{1}(0)=0$,  and its largest speed,
 obtained at $t\in[10,15]$ and again at $t\in[25,30]$,
 is $20$m/s. At the point of largest curvature, when $z_{2}$ attains its maximum (see Figure 6), its speed is $8.66$m/s.  The four vehicles start at rest at the point $r_{1}(0)$, and the initial condition of their controller is $u(0)=(2,0)^{\top}$.

Figures 6-9 present simulation results with the controller
defined by  Eq. \eqref{eq:control_short}.
Figure 6 depicts the target and agent-trajectories  from left to right in the $(z_{1},z_{2})$ plane. Both coordinates $z_{1}$ and $z_{2}$ are of the same scale
thereby indicating  quite large curvature of the target trajectory at the point of maximum $z_{2}$.
Figure  8 shows the graphs of the lateral (normal)  errors of the vehicles' centers of gravity   from the target trajectory $\{r_{1}(t)\}$, and we note that  the relatively large error-spurts correspond to the larger  curvatures
indicated  in Figure 6. Furthermore, as expected, the   errors of later vehicles in the platoon tend to be larger than those of earlier ones.  The maximum lateral error,
obtained for $A_{4}$, is about 38 cm.

Graphs of an approximate measure of the  inter-agent distances vs. time  are shown in Figure 9. We have to use an approximate (not exact) distance for the following
reason: The objective of the control law is to drive the vehicles to  the path $\{r_{1}(t)\}$ where they maintain an inter-agent distance of $10$m. The term ``distance'' between two consecutive vehicles means the arclength between them, which is well defined as long as both vehicles are on the path, but not well defined when one or both of them are off the path.  Therefore we display, in Figure 9,  the approximate measure of distance between two vehicles defined as the sum of the Euclidean distance of each vehicle to the nearest-point to it on the path, and the arclength between these two nearest points.\footnote{The nearest point is assumed to be unique.}  The justification for this measure of distance is that the control algorithm drives the vehicles towards the path, where this measure coincides with the arclength. In fact, Figure 8 shows that the vehicles converge to the path $\{r_{1}\}$ except for at points of large curvature, and Figure 9 displays a convergence of the corresponding measure of distance towards 10m except at such points.

Finally, the longitudinal accelerations of the vehicles are depicted in  Figure 7. Although they may make for an  uncomfortable ride,  they closely track  the acceleration of the  target path $\{r_{1}(t)\}$, with a notable  deviation corresponding  to its region of largest curvature.

\begin{figure}	
	\centering
		\includegraphics[width=0.75\linewidth]{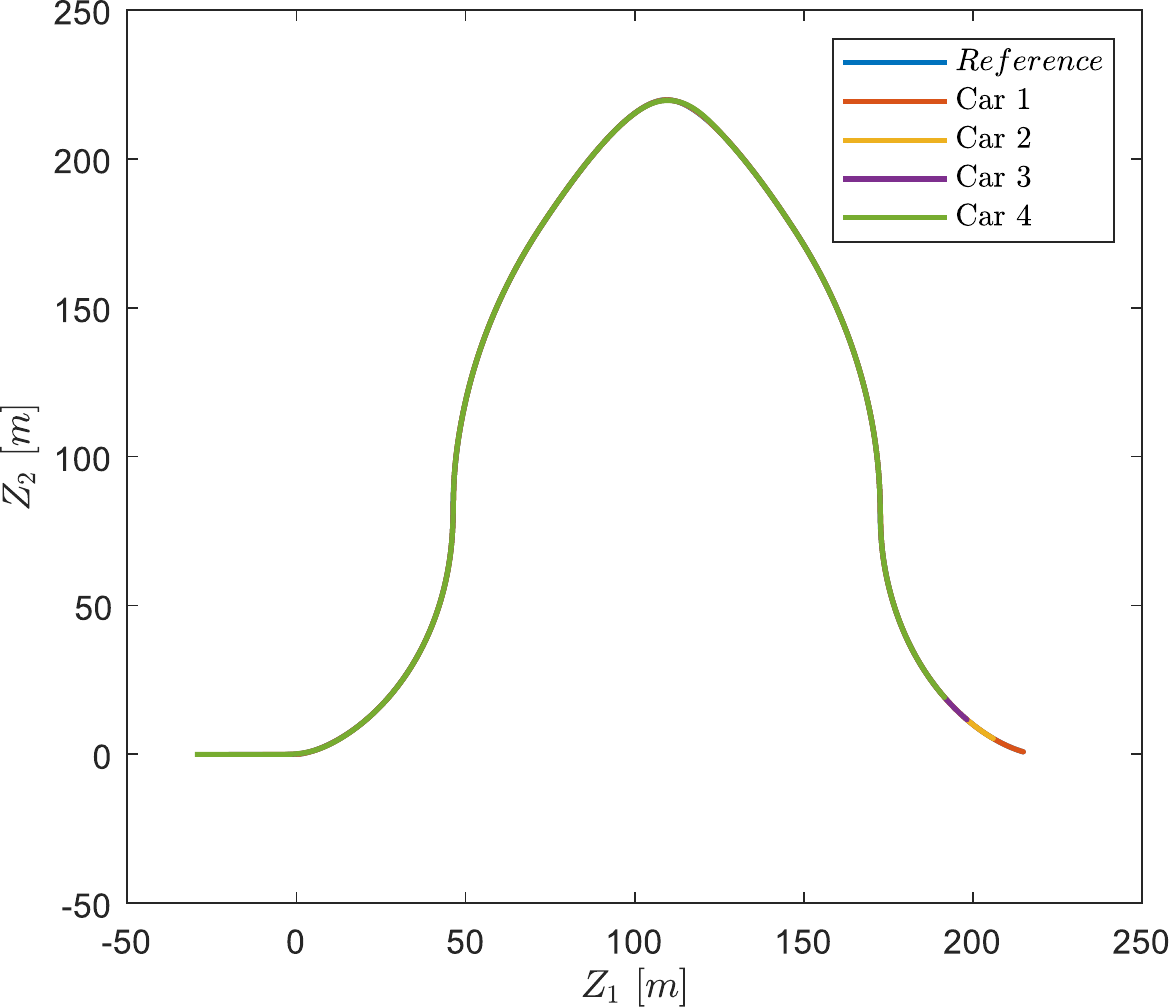}
	\caption{{\small Platoon: target trajectory in the $z$-plane}}
\end{figure}

\begin{figure}
	\centering
	\includegraphics[width=0.75\linewidth]{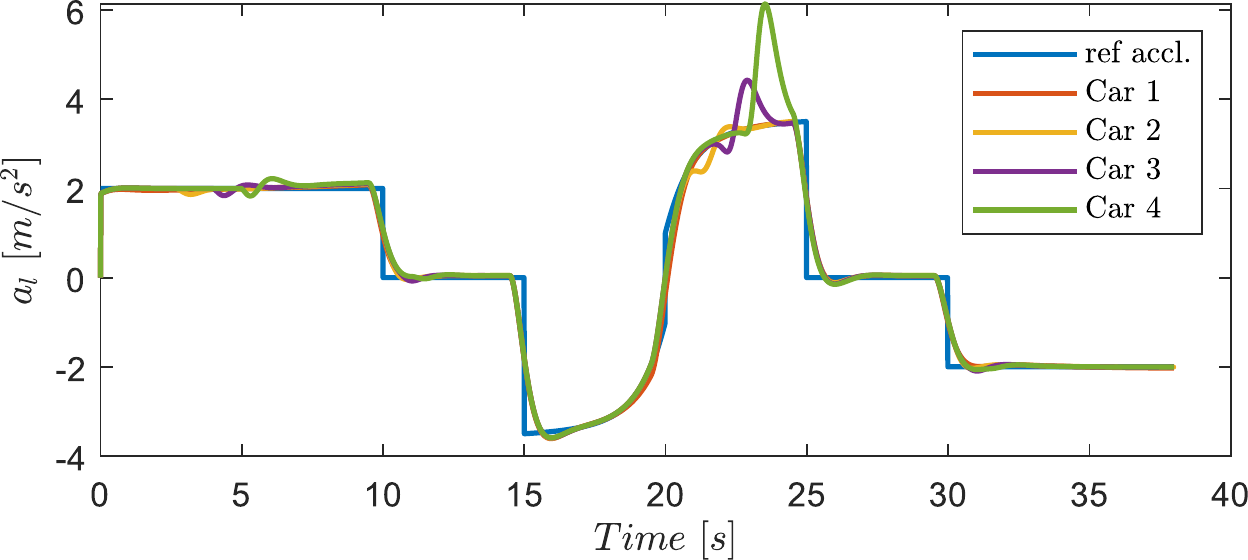}
	\caption{{\small Platoon: Reference-path and vehicle  accelerations}}
\end{figure}

 \begin{figure}
	\centering
	\includegraphics[width=0.75\linewidth]{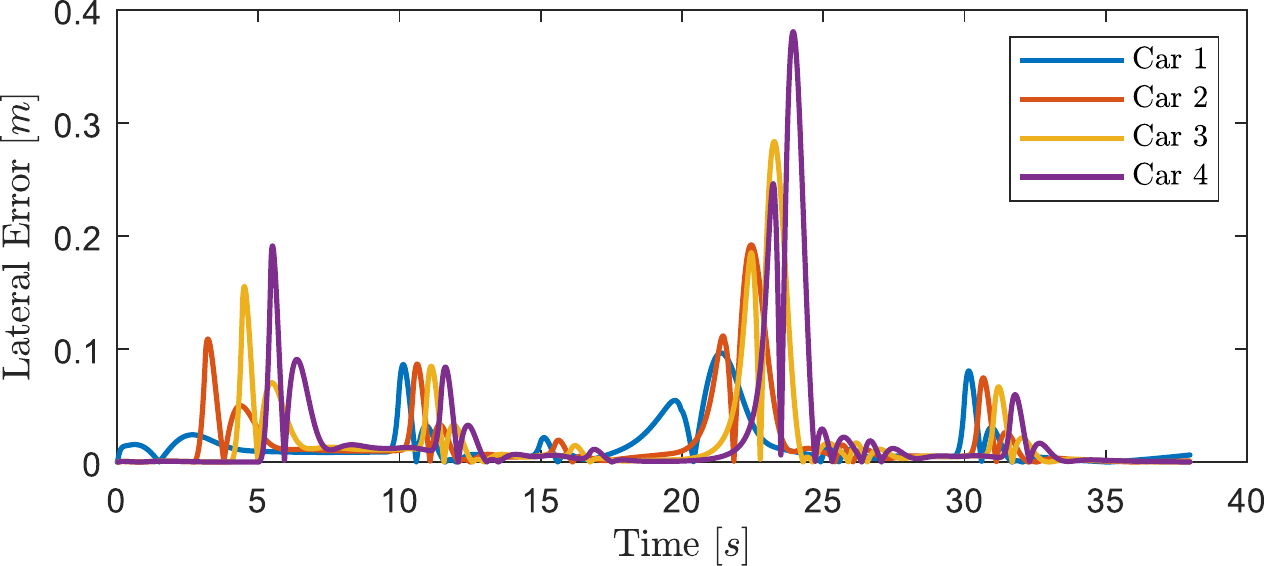}
	\caption{{\small Platoon: lateral errors vs. time}}
\end{figure}

  \begin{figure}  
	\centering
	\includegraphics[width=0.75\linewidth]{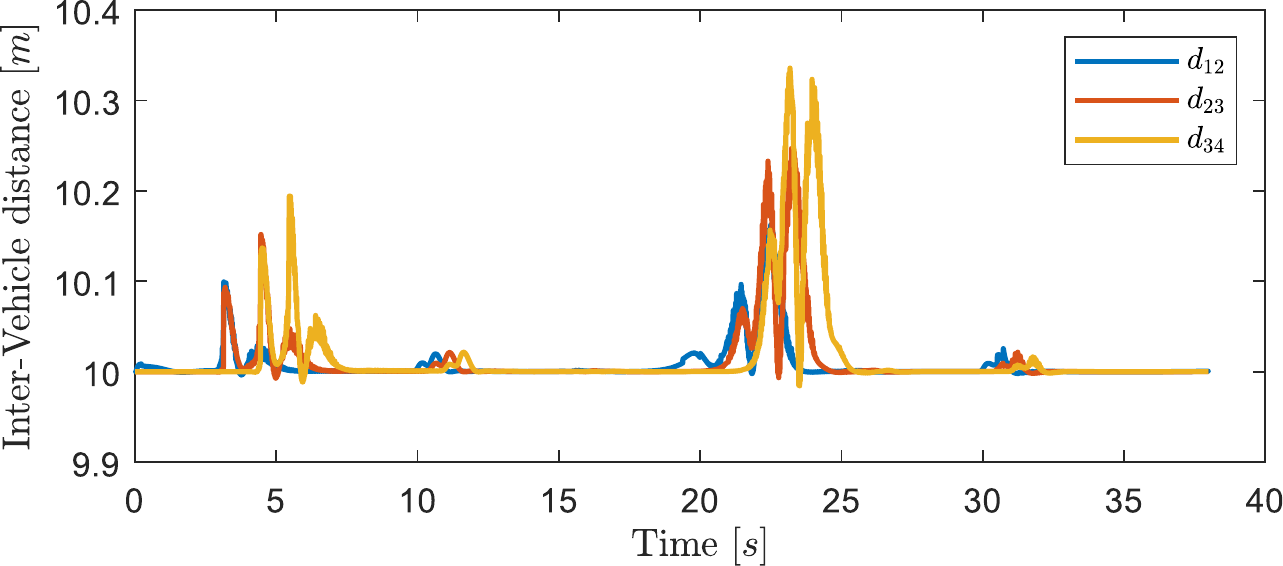}
	\caption{{\small Platoon: approximate inter-agent distances}}
\end{figure}

\begin{figure}  
    \centering
    \includegraphics[width=0.85\linewidth]{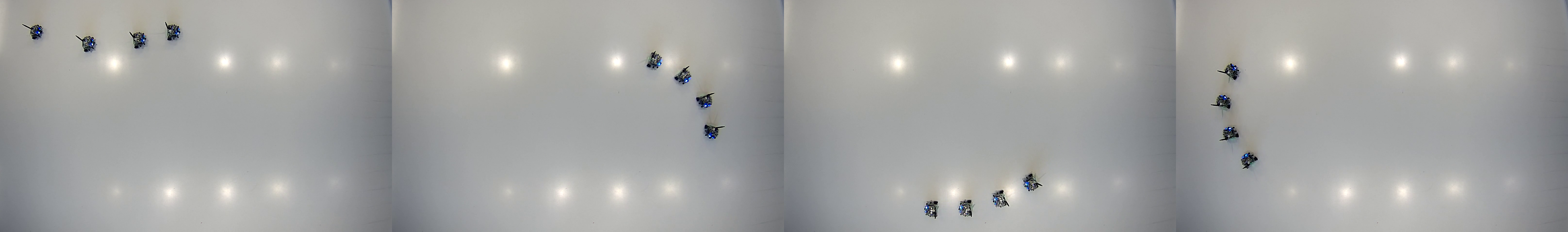}
    \caption{{Experiment: Stills of the robots' positions around the curve}}
\end{figure}

\section{Experimental Results}

This section describes results of laboratory experiments in which a platoon of four mobile robots (agents)  attempts to maintain a given inter-agent   distance. The  present system is different from the one considered in Subsection 5.2 in several ways including the following three:  (i) The experimental setting is a laboratory vs.   simulation, (ii) the vehicles' dynamic equations follow a unicycle model vs. a bicycle model,
and (iii) $A_{i}$, $i=2,3,4$, only have to maintain the given inter-agent distance from $A_{i-1}$ but not follow its trajectory.

The platoon consists of four agents denoted by  $A_{i}$, $i=1,\ldots,4$, according to their order. The lead agent, $A_{1}$, is assigned its planar target trajectory, $\{r_{1}(t)\}$,  by an exogenous source,
and for every $i=2,3,4$,  $A_{i}$ aims at keeping a given
Euclidean distance from $A_{i-1}$.

The experiments were conducted in the Robotarium, a
remotely-accessible testing facility for motion control of robotic systems located at the Georgia Tech campus \cite{Pickem17}. The vehicles in the Robotarium are  differential-drive robots, approximately 15cm in diameter, which were designed and assembled in-house. Their motion is  modelled by unicycle dynamics having the following form,
\begin{equation}
\begin{pmatrix}
\dot{z}_{1}(t)\\\dot{z}_{2}(t)\\\dot\psi(t)
\end{pmatrix} = \begin{pmatrix}
\cos\psi(t) & 0\\ \sin\psi(t) & 0 \\ 0 & 1
\end{pmatrix}\begin{pmatrix}
v(t) \\ \omega(t)
\end{pmatrix},
\label{eq:unicycle}
\end{equation}
where $z:=(z_{1},z_{2})^{\top}\in R^2$ is the  center of gravity of a robot and $\psi$ is  its  heading. Eq. (75) is a state-space representation of a vehicle with the state variable  $x:=(z_{1},z_{2},\psi)^{\top}$ and a control input  $u:=(v,\omega)^{\top}$, where  $v$  and $\omega$ are  its longitudinal velocity and  angular velocity, respectively. The output of the system is $y(t):=z(t)=(z_{1}(t),z_{2}(t))^\top$.

Fix a prediction horizon $T>0$. A direct integration of
Eq. (17),
together with (18), result in the following closed-form
for the output predictor $\hat{y}(t+T):=g(x(t),u(t))$,
\begin{equation}
    g(x(t),u(t))=
    \begin{pmatrix}\begin{array}{c}z_{1}(t)\\z_{2}(t)\end{array}\end{pmatrix}
    +\frac{v(t)}{\omega(t)}\begin{pmatrix}\begin{array}{c}
    \sin\big(\psi(t)+\omega(t)T\big)-\sin\big(\psi(t)\big)\\
    -\cos\big(\psi(t)+\omega T\big)+\cos\big(\psi(t)\big)
    \end{array}
  \end{pmatrix};
  \label{eq:predictor}
\end{equation}
 if $\omega(t)=0$,  L'Hopital's rule  yields
 \begin{equation}
    \label{eq:lhopital}
     g(x(t),u(t))=
     \begin{pmatrix}\begin{array}{c}z_{1}(t)\\z_{2}(t)\end{array}\end{pmatrix}+v(t)T\begin{pmatrix}\begin{array}{c}
     \cos(\psi(t))\\\sin(\psi(t))
     \end{array}
     \end{pmatrix}.
 \end{equation}
 The controller uses this functional closed form and does not  resort to numerical integration of (17).

 The future target-point  $r_{i}(t+T)$ is defined for  the agent $A_{i}$ according to the following heuristic.
  For $i=1$,  $\{r_{1}(t)\}$ is assumed to be known in advance and hence $r_{1}(t+T)$ can be used in the computations of $A_{1}$  at time $t$.
For $i=2,3,4$, the definitions and computations of  $r_{i}(t+T)$ are recursive, as follows.  At time $t$, let   $\ell_{i}$ denote the directional line from $\hat{y}_{i-1}(t+T)$  towards $y_{i}(t)$, namely the line connecting the predicted position of $A_{i-1}$ towards the current position of $A_{i}$. Then we define   $r_{i}(t+T)$  as  the point on $\ell_{i}$ of distance $d$ m from $\hat{y}_{i-1}(t+T)$.   This procedure is justified by the observation that  if $A_{1}$ moves in a straight line,  then subsequent agents will converge to that line behind each other at the target distance $d$.

We conducted experiments with the controller defined by Eq.   (38), $\alpha=45$ and $T=0.25$s.
The exogenous target curve, $\{r_{1}(t)\}$, is an ellipse defined by the equation
$r_1(t) =
            \big(1.1 \sin(0.06 t),0.7 \cos(0.06 t)\big)^{\top},$
and the target inter-robot distance is $d=0.25$m.
The results are depicted in Figures 10-12.
Figure 10 shows stills captured during the experiment. In the leftmost image the robots are initialized, and in subsequent images of their positions are shown; the first robot moves along the closed curve defined by $\{r_{1}(t)\}$ while the remaining robots converge to the target inter-robot distances.  Figure 11 depicts  the inter-robot distances vs. $t$; note convergence towards the target distance of $0.25$m. Finally,  Figure 12 depicts the graph of the tracking error $\|y_{i}(t)-r_{i}(t)\|$ versus time, and we discern rapid convergence towards 0 for all four robots. An additional view of the control-algorithm's performance can  be seen in the video clip contained in \cite{Buckley19}.

\begin{figure}
    \centering
    \includegraphics[width=0.75\linewidth]{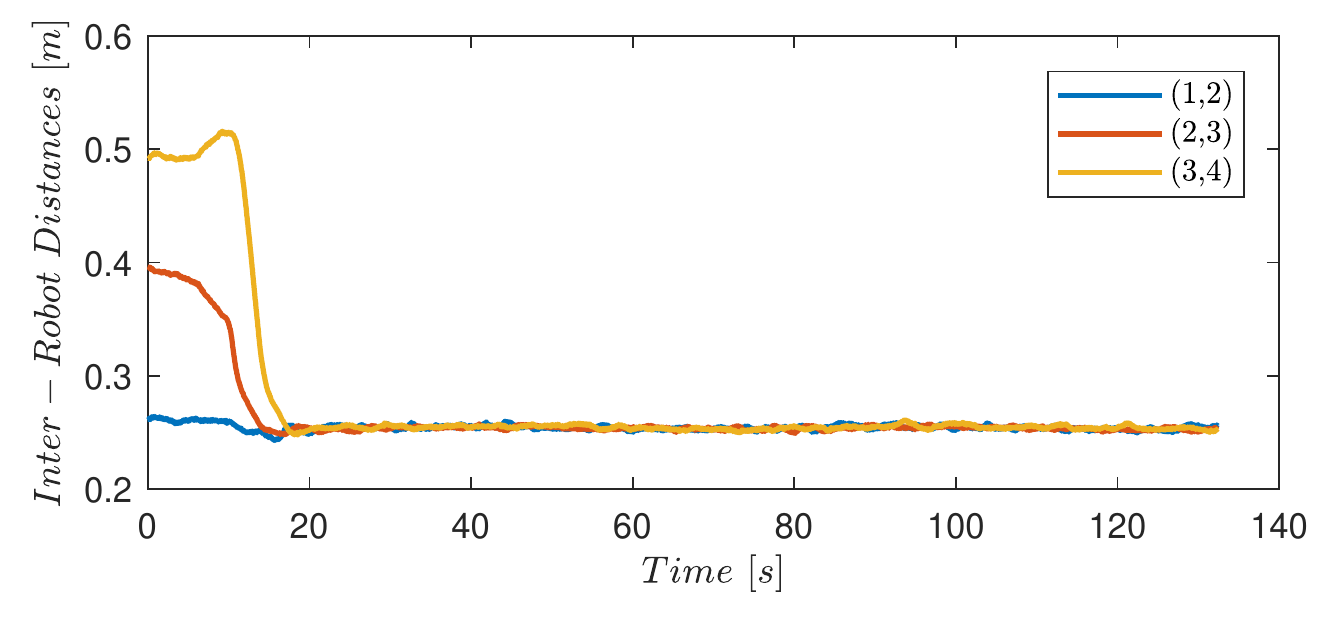}
    \caption{{\small Experiment: inter-robot distances vs. time }}
    \label{fig:experiment-distance}
\end{figure}

\begin{figure}
    \centering
    \includegraphics[width=0.75\linewidth]{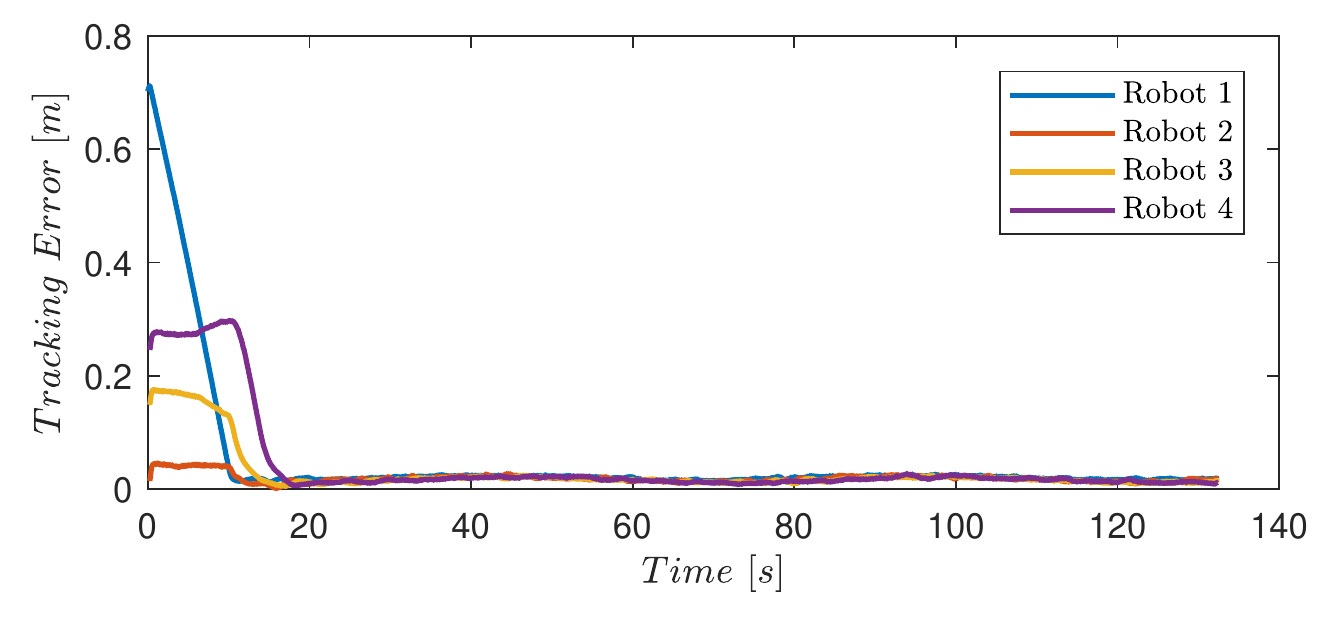}
    \caption{{\small Experiment: tracking error vs. time}}
    \label{fig:experiment-tracking}
\end{figure}

\section{Conclusions}
This paper presents a tracking-control technique based on a fluid-flow version of the Newton-Raphson method, output prediction and controller speedup. The controller is simple to compute  and may have large, even global stability domains of attraction.  A stability analysis is carried out for linear systems, while examples of nonlinear systems are tested by simulation and lab experiments.

Current investigations concern theoretical and practical problems. On the theoretical side, the most pressing challenge is to derive sufficient conditions for the $\alpha$ stability of  closed-loop systems in a general setting  of nonlinear dynamical systems. Practical considerations include the testing of the control technique on applications in mobile robotics and autonomous vehicles. Of a particular interest is to derive alternative output-prediction techniques to the one presented in this paper, and what comes to mind are methods that are based on learning and neural nets.

\section{Appendix}
This section provides  proofs of various assertions made throughout the paper.

{\it Proof of Lemma 4.3.} The characteristic polynomial of $\Phi_{\alpha}$, denoted by $P_{\alpha}(s)$, is a two-dimensional polynomial in $(\alpha,s)$. Denote its respective degrees in $\alpha$ and
$s$ by $deg_{\alpha}(P)$ and  $deg_{s}(P)$,
and define its total degree, denoted by $deg(P)$, as the degree
(in $q$) of the polynomial $P_{q}(q)$.

By Eq. (49), we observe that every element (entry) in the matrix
$sI-\Phi_{\alpha}$ has a total degree of zero or 1; for example, the $(n+m)\times(n+m)$ element is $s-\alpha\phi_{n+m,n+m}$ which contains both $\alpha$ and $s$ but not the product of the two. Since the determinant of a matrix consists of the linear combination of products of elements one from each row, we have that
\begin{equation}
    deg(P)=n+m.
\end{equation}
Furthermore, $deg_{\alpha}(P)=m$, since only the last $m$ rows of $\Phi_{\alpha}$ contain the term $\alpha$. Therefore, Eq. (50) is in force for some polynomials $P_{m-i}(s)$, $i=0,\ldots,m$, and by 78), the degree (in $s$) of $P_{m-i}$ must not exceed $n+(m-i)$. This completes the proof. \\ \\
{\it Proof of Lemma 4.6.}
We have that $P_{\alpha}(s(\alpha))=0$ $\forall~\alpha\in[0,\infty)$. Therefore, and by Eq. (50),
 \begin{equation}
 \sum_{i=0}^{m}\alpha^i P_{m-i}(s(\alpha))=0.
 \end{equation}
 Dividing the latter equation by $\alpha^m$, we obtain that
 \begin{equation}
 \sum_{i=0}^{m-1}\alpha^{i-m}P_{m-i}(s(\alpha))+P_{0}(s(\alpha))=0.
 \end{equation}
 Since $\{s(\alpha)\}$ is bounded, the sum-term in the RHS of Eq. (80)  goes to 0 as $\alpha\rightarrow\infty$. Therefore,
 taking $\alpha\rightarrow\infty$ in (80), we have that
 \[
 \lim_{\alpha\rightarrow\infty}P_{0}(s(\alpha))=0.
 \]
 Since $\{s(\alpha)\}$ is bounded, it has at least one limit (accumulation) point;
 and by the latter equation, such a limit point must be a root of $P_0(s)$.
 Since $P_0(s)$ has a finite number
 ($n$) of  roots, the limit
 $\lim s(\alpha)$ (as
 $\alpha\rightarrow\infty$) exists  and it is a root of $P_0(s)$.
\hfill $\Box$\\ \\
{\it Proof of the left inequality of Eq. (59).}
 We argue by contradiction. If the left inequality in (59) is not satisfied, there exists an unbounded
 set $A_{2}\subset A$ such that, as $\alpha\rightarrow\infty;~\alpha\in A_{2}$,
 \begin{equation}
 \frac{|s(\alpha)|}{\alpha}\rightarrow 0.
 \end{equation}
 By (50), for every $\alpha\in A_{2}$,
 \[
 \sum_{i=0}^{m}\alpha^i P_{m-i}(s(\alpha))=0.
 \]
 Divide this equation by $\alpha^m s(\alpha)^n$ to obtain, $\forall~\alpha\in A_{2}$,
 \begin{equation}
 \sum_{i=0}^{m-1}\frac{P_{m-i}(s(\alpha))}{\alpha^{m-i}s(\alpha)^{n}}+\frac{P_{0}(s(\alpha))}{s(\alpha)^n}
 =~\sum_{i=0}^{m-1}\frac{s(\alpha)^{m-i}}{\alpha^{m-i}}
 \cdot\frac{P_{m-i}(s(\alpha))}{s(\alpha)^{m+n-i}}+\frac{P_{0}(s(\alpha))}{s(\alpha)^n}~=~0.
 \end{equation}
 By Eq. (81), and since ${\rm deg}(P_{m-i})=n+m-i$ for all $i=0,\ldots,m-1$; as $\alpha\rightarrow\infty,~\alpha\in A_{2}$,
 \[
\sum_{i=0}^{m-1}\frac{s(\alpha)^{m-i}}{\alpha^{m-i}}\cdot\frac{P_{m-i}(s(\alpha))}{s(\alpha)^{m+n-i}}~\rightarrow~0.
 \]
 Furthermore, since ${\rm deg}(P_{0})=n$,
 \[
 \lim_{\alpha\in A_2; \alpha\rightarrow\infty}\frac{P_{0}(s(\alpha))}{s(\alpha)^n}\neq 0.
 \]
 This contradicts Eq. (82) and hence completes the proof.   \hfill$\Box$\\ \\
 {\it Proof of Lemma~\ref{le:lemma4}.}
  By Eq. (52), for every $\alpha\geq 0$,
 \begin{equation}
 P_{\alpha}(s(\alpha))=\sum_{i=0}^{m}\alpha^{i}\sum_{j=0}^{n+m-i}a_{m-i,j}s(\alpha)^{j}=0.
 \end{equation}
Fix $\ell\in\{0,\ldots,m\}$ and $\nu\in\{0,\ldots,n+m-\ell -1\}$. Taking derivatives in (83) with respect to
$a_{m-\ell,\nu}$ we obtain,
\begin{equation}
\sum_{i=0}^{m}\sum_{j=0}^{n+m-i}\alpha^{i}a_{m-i,j}\cdot js(\alpha)^{j-1}\frac{\partial s(\alpha)}{\partial a_{m-\ell,\nu}}+\alpha^{\ell}s(\alpha)^{\nu}=0,
\end{equation}
hence
\begin{equation}
\frac{\partial s(\alpha)}{\partial a_{m-\ell,\nu}}=-\frac{\alpha^{\ell}s(\alpha)^{\nu}}{\sum_{i=0}^{m}
\sum_{j=0}^{n+m-i}\alpha^{i}a_{m-i,j}\cdot j s(\alpha)^{j-1}}.
\end{equation}
Now both numerator and denominator in Eq. (85) are comprised of two-dimensional
polynomials in $\alpha$ and $s=s(\alpha)$. Their total degrees are
$n+m-1$ for the denominator, and $\ell+\nu$ for the numerator. But
$\nu\leq n+m-\ell-1$ by assumption, hence $\ell+\nu\leq n+m-1$, implying that
the total degree of the numerator is less or equal to that of the denominator.
This, together with Lemma~\ref{le:lemma3}, imply that Eq. (62) and hence the lemma'a assertion. \hfill $\Box$

\bibliographystyle{IEEEtran}

\bibliography{Paper}

\end{document}